\documentclass{imsart}
\RequirePackage[OT1]{fontenc}
\RequirePackage{amsthm,amsmath}
\RequirePackage[colorlinks,citecolor=blue,urlcolor=blue]{hyperref}

\usepackage{amsmath}
\usepackage{amssymb}
\usepackage{subfigure}
\usepackage{graphicx}
\usepackage{epstopdf}
\usepackage{bm}
\usepackage{mathrsfs}
\usepackage{enumerate}
\usepackage{multirow}
\usepackage{array}
\usepackage{color}
\usepackage{threeparttable}
\usepackage[round]{natbib}

\newcommand{\cA}{{\mathcal A}}

\renewcommand{\theenumi}{\Alph{enumi}}

\makeatletter
\renewcommand{\p@enumii}{\theenumi.}
\makeatother
\def\namedlabel#1#2{\begingroup
    #2%
    \def\@currentlabel{#2}%
    \phantomsection\label{#1}\endgroup
}
\makeatother

\def\@enum@{\list{\csname label\@enumctr\endcsname}%
          {\usecounter{\@enumctr}\def\makelabel##1{\hss\llap{##1}}
          \itemsep=2pt\parsep=0pt\topsep=3pt plus 1pt minus 1 pt}}
\usepackage{paralist}

\newtheorem{thm}{Theorem}[section]
 \newtheorem{cor}{Corollary}[section]
 \newtheorem{lem}{Lemma}[section]
 \newtheorem{prop}{Proposition}[section]
 \theoremstyle{definition}
 \newtheorem{defn}{Definition}[section]
 \newtheorem{rmk}{Remark}[section]

\newcommand{\bA}{{\bf A}}

\newcommand{\bR}{{\bf R}}
\newcommand{\bM}{{\bf M}}
\newcommand{\bX}{{\bf X}}
\newcommand{\bY}{{\bf Y}}
\newcommand{\bZ}{{\bf Z}}
\newcommand{\bW}{{\bf W}}
\newcommand{\bV}{{\bf V}}
\newcommand{\bU}{{\bf U}}

\newcommand{\bS}{{\bf S}}

\newcommand{\bI}{{\bf I}_p}

\newcommand{\bP}{\widetilde{\bf U}}

\newcommand{\bp}{\widetilde{\bf u}}

\newcommand{\cC}{{\mathcal C}}

\newcommand{\mT}{{\mathrm{T}}}
\newcommand{\mbA}{{\tilde{\mathcal{\bm A}}}}  
\newcommand{\mbB}{{\mathcal{\bm B}}}

\newcommand{\bv}{{\bf v}}
\newcommand{\bw}{{\bf w}}

\newcommand{\bXi}{\boldsymbol{\Xi}}

\newcommand{\bmu}{\boldsymbol{\mu}}
\newcommand{\Var}{\textrm{Var}}
\newcommand{\RCV}{\textrm{RCV}}
\newcommand{\ICV}{\textrm{ICV}}
\newcommand{\TVA}{\textrm{TVA}}
\newcommand{\MNS}{\textrm{MNS}}
\newcommand{\APA}{\textrm{APA}}
\newcommand{\ANS}{\textrm{ANS}}
\newcommand{\diag}{\textrm{diag}}

\newcommand{\bbSigma}{\breve{\bSigma} }
\newcommand{\wbSigma}{\widetilde{\bSigma}}

\newcommand{\bSigma}{\boldsymbol{\Sigma}}
\newcommand{\bLambda}{\boldsymbol{\Lambda}}
\newcommand{\bTheta}{\boldsymbol{\Theta}}
\def\hn{N}
\def\hv{M}
\def\Var{\begin{rm} Var \end{rm}}

\def\no {\noindent}

\def \proofend{\hspace{\fill}$\blacksquare$ \hbox{} \medskip}
\def\widebar{\accentset{{\cc@style\underline{\mskip10mu}}}}
\def\widebar{\overline}

\def\wwbY{\widetilde{\overline{\bf{Y}}}}
\def\twbY{\widetilde{\overline{\bf{Y}}}}
\def\twbX{\widetilde{\overline{\bf{X}}}}
\def\twbe{\widetilde{\overline{\bm{\varepsilon}}}}
\def\twbV{\widetilde{\overline{\bf{V}}}}
\def\wbX{\overline{\bf{X}}}
\def\wbV{\overline{\bf{V}}}

\def\tw{\widetilde{w}}
\def\ttw{\widetilde{\widetilde{w}}}

\def\ttbSigma{\widetilde{\widetilde{\bm{\Xi}}}}

\def\hbSigma{\widetilde{\bm{\Xi}}}

\def\tbZ{\widetilde{\bf{Z}}}
\def\tbV{\widetilde{\bf{V}}}

\def\tr{\begin{rm} tr\end{rm}}

\def\ol{\overline}

\def\wt{\widetilde}
\def\wh{\widehat}
\def\indic{{\mathbf{I}}}

\startlocaldefs
\numberwithin{equation}{section}
\theoremstyle{plain}

\endlocaldefs
\allowdisplaybreaks[4]

\begin{document}
\pagenumbering{Alph}
\begin{titlepage}
\maketitle
\thispagestyle{empty}
\end{titlepage}
\pagenumbering{arabic}

\thispagestyle{empty}

\begin{frontmatter}
\title{On the estimation of high-dimensional integrated covariance matrices based on high-frequency data with multiple transactions}
\runtitle{Estimation of ICV based on multiple noisy obs}

\begin{aug}
\author{\fnms{Moming} \snm{Wang}\thanksref{t1}\ead[label=e1]{wangmoming0902@163.com}},
\author{\fnms{Ningning} \snm{Xia}\thanksref{t2}\ead[label=e2]{xia.ningning@mail.shufe.edu.cn}}
\and
\author{\fnms{Yong} \snm{Zhou}\thanksref{t3}
\ead[label=e3]{yzhou@amss.ac.cn}}

\thankstext{t1}{Research partially supported by the NSFC 11871322.}
\thankstext{t2}{Research partially supported by the NSFC 11871322 and 11571218, IRTSHUFE.}
\thankstext{t3}{Research partially supported by the State Key Program of NSFC 71331006, the State Key Program in the Major Research Plan of NSFC 91546202.}
\runauthor{M. Wang, N. Xia and Y. Zhou}

\affiliation{Shanghai University of Finance and Economics\thanksref{t1,t2},
East China Normal University\thanksref{t3} }

\address{School of Statistics and Management, Shanghai University of Finance and Economics,\\
777 Guo Ding Road, Shanghai, 200433, China\\
\printead{e1}\\
}

\address{School of Statistics and Management, Shanghai, Key Laboratory of Financial Information Technology,
Shanghai University of Finance and Economics, \\
777 Guo Ding Road, Shanghai, 200433, China\\
\printead{e2}\\
}

\address{Faculty of Economics and Management, East China Normal University,\\ Shanghai, 200241, China
Academy of Mathematics and Systems Science, Chinese Academy of Science, Beijing, 100190, China\\
\printead{e3}}
\end{aug}

\begin{abstract} 
High-frequency data in financial markets often include multiple transactions at each recording time due to the mechanism of recording.
Using random matrix theory, this paper considers the estimation of integrated covariance (ICV) matrices of high-dimensional diffusion processes based on multiple high-frequency observations.
We start by studying the estimator, the time-variation adjusted realized covariance (TVA) matrix proposed by \cite{ZhengandLi2011}, without microstructure noise.
We show that in the high-dimensional case, for a class $\mathcal{C}$ of diffusion processes, the limiting spectral distribution (LSD) of the averaged TVA depends not only on that of the ICV but also on the number of multiple transactions at each recording time.
However, in practice, observed prices are always contaminated by market microstructure noise.
Thus, we also study the limiting behavior of pre-averaging averaged TVA matrices based on noisy multiple observations.
We show that for processes in class $\mathcal{C}$, the pre-averaging averaged TVA has two desirable properties: it eliminates the effects of microstructure noise and multiple transactions, and its LSD depends solely on that of the ICV matrix.
Further, three types of nonlinear shrinkage estimators of the ICV matrix are proposed based on high-frequency noisy multiple observations.
Simulation studies support our theoretical results and demonstrate the finite sample performance of the  proposed estimators.
Finally, high-frequency portfolio strategies are  evaluated under these estimators in real data analysis.

\end{abstract}



\end{frontmatter}

\section{Introduction}\label{section: introduction}

The covariance structure of capital markets is of great interest to investors and researchers, as it has a critical role in financial problems such as pricing and investment.
Suppose that we have $p$ stocks whose latent log price at time $t$ is denoted by $\bX_t = (X_t^{(1)}, \cdots, X_t^{(p)})^\mathrm{T}$, where $T$ denotes the transpose.
The following diffusion processes are commonly used to model financial asset price processes:
\begin{equation}
d\bX_t = \bmu_tdt + \bTheta_td\bW_t, \hspace{1cm} t\in [0,1],
\label{eq: Model}
\end{equation}
where $(\bmu_t) = (\mu_t^{(1)},\cdots,\mu_t^{(p)})^\mathrm{T}$ is a $p$-dimensional drift process, $(\bTheta_t)$ is a $p\times p$ matrix, the so-called covolatility process at time $t$, and $(\bW_t)$ is a $p$-dimensional standard Brownian motion.
The interval $[0,1]$ represents the time period of interest, such as one trading day.
The integrated covariance (ICV) matrix given by
\begin{equation}
\ICV :=\int_0^1 \bTheta_t\bTheta_t^\mathrm{T} dt
\label{eq: definition of ICV}
\end{equation}
is of fundamental importance in risk management and portfolio allocation for high-frequency financial data.
In the one-dimensional case, the ICV matrix reduces to the integrated volatility.
The estimation of the ICV matrix is a fundamental problem in financial applications.

The classical estimator of the ICV matrix is the so-called realized covariance (RCV) matrix, which is defined as follows. Suppose that we can observe the log price processes $\{X_{t_i}^{(j)}\}_{j=1}^p$ synchronously at recording time $t_i = i/n$, $i=0,1,\cdots,n$, during one trading day. Then the RCV matrix is defined as
\[\RCV := \sum_{i=1}^n\Delta\bX_i(\Delta\bX_i)^{\mathrm{T}},\]
where
\[
\Delta\bX_i = \begin{pmatrix}
\Delta X_i^{(1)}\\
\vdots \\
\Delta X_i^{(p)}\\
\end{pmatrix}
=
\begin{pmatrix}
 X_{t_i}^{(1)} -  X_{t_{i-1}}^{(1)}\\
\vdots\\
 X_{t_i}^{(p)} -  X_{t_{i-1}}^{(p)}
\end{pmatrix}.
\]
In the case where dimension $p$ is fixed and the number of observations $n$ goes to infinity, the consistency and central limit theorem for the RCV have been well studied, for example, by \cite{AndersenandBollerslev1998}, \cite{Andersenetal2001}, \cite{BNS2002}, and \cite{JacodandProtter1998}.
However, the classical estimator RCV matrix is no longer consistent in high-dimensional settings, where the number of stocks $p$ increases with the observation frequency $n$ at the same rate (see \cite{WangandZou2010}, \cite{Fanetal2012}, \cite{Taoetal2011}, \cite{ZhengandLi2011}, and \cite{XiaandZheng2018}, for instance).
Notably, when the sparsity assumptions on the high-dimensional ICV matrix are not satisfied, a key assumption can be used to solve the problem of estimation of the ICV matrix, that is, the class $\mathcal{C}$ condition.

\begin{defn}
Suppose that $(\bX_t)$ is a $p$-dimensional process satisfying (\ref{eq: Model}) and $\bTheta_t$ is c\`{a}dl\`{a}g. Then $(\bX_t)$ belongs to class $\mathcal{C}$ if, almost surely, there exist $\gamma_t\in D([0,1]; \mathbb{R})$ and $\bLambda$ a $p\times p$ matrix satisfying $\begin{rm}tr\end{rm}(\bLambda\bLambda^\mathrm{T}) = p$ such that
\begin{equation}
\bm\Theta_t = \gamma_t \bm \Lambda,
\label{eq: definition of Theta}
\end{equation}
where $D([0,1]; \mathbb{R})$ is the space of c\`{a}dl\`{a}g functions from $[0,1]$ to $\mathbb{R}$.
\end{defn}
It has been shown that the log price process $(\bX_t)$ belongs to class $\cC$ in Proposition 4 of \cite{ZhengandLi2011}, with the reasonable assumption that the correlation structure of $(\bX_t)$ does not change over a short time period.
Moreover, the class $\cC$ process accommodates both stochastic volatility and the leverage effect, which are common features of financial data.
Thus, for class $\cC$ processes, the limiting behavior of the RCV matrix has been well studied, by \cite{ZhengandLi2011}, \cite{Lametal2017}, \cite{XiaandZheng2018}, and others.
However, the RCV matrix cannot be used to make robust inferences regarding the spectrum of the ICV matrix, owing to the difficulty of estimating the unknown and time-volatile covolatility process $(\gamma_t)$.
Hence, an alternative estimator has been proposed, the time-variation adjusted realized covariance (TVA) matrix, which is defined as
\[
\TVA:= \frac{\tr(\RCV)}{n}\sum_{i=1}^n \frac{\Delta\bX_i(\Delta\bX_i)^{\mT}}{|\Delta\bX_i|^2}.\]
The basic idea is to eliminate the influence of the (unknown) covolatility process $(\gamma_t)$ from the limiting spectral relationship between the  ICV matrix and its estimator TVA matrix.
Although the TVA matrix remains inconsistent, using random matrix theory,
Theorem 2 of \cite{ZhengandLi2011} shows that
its limiting spectral distribution depends solely on the ICV matrix through the Mar\v{c}enko--Pastur equation \eqref{eq: m-p equation} (see Theorem \ref{thm: theorem for the TVAPA}).  Therefore, the spectrum of the ICV can be recovered using existing algorithms, such as those developed by \cite{Karoui2008}, \cite{Mastre2008}, and \cite{Baietal2010}.

However, these approaches have only been applied to datasets that contain exactly one transaction during one time stamp.
Specifically, the estimators,  RCV and TVA, are
established based on the assumption that there is only one observation $\bX_{t_i}$ at each recording time $t_i$.
However, this is not necessarily the case in practice.
In high-frequency financial markets, owing to heavy market trading and the limitations of the recording mechanism, multiple transactions often occur at each recording time.
\begin{table}[!htb]
\small
\centering
\setlength{\abovecaptionskip}{0pt}
\setlength{\belowcaptionskip}{10pt}
\caption{Numbers of transactions (NT) and different prices (NDC) of Apple, Cisco  Systems and Microsoft in the first 30 seconds on November 4, 2016 }
\setlength{\tabcolsep}{5 mm}{
\begin{tabular}{lllllll}
  \hline
  \hline
 &\multicolumn{2}{c}{Apple}	&\multicolumn{2}{c}{Cisco Systems}	&\multicolumn{2}{c}{Microsoft}\\

	Time &NT	& NDC	&NT	& NDC	&NT	& NDC\\
\hline
    9:30:01	&492	&42	&166	&13	&217	&30\\
    9:30:02	&357	&28	&232	&37	&186	&22\\
    9:30:03	&442	&31	&120	&14	&150	&23\\
    9:30:04	&270	&54	&46	    &15	&134	&18\\
    9:30:05	&288	&30	&69	    &13	&86	   &13\\
    9:30:06	&163	&28	&100	&26	&149	&17\\
    9:30:07	&230	&47	&212	&9	&239	&24\\
    9:30:08	&289	&54	&104	&7	&91	   &20\\
    9:30:09	&80	    &31	&82	    &10	&145	&16\\
    9:30:10	&94	    &33	&77	    &10	&122	&14\\
    9:30:11	&129	&35	&149	&14	&212	&14\\
    9:30:12	&222	&26	&82	    &17	&107	&8\\
    9:30:13	&126	&33	&9	    &8	&268	&11\\
    9:30:14	&106	&19	&16	    &6	&102	&10\\
    9:30:15	&131	&34	&65	    &13	&130	&12\\
    9:30:16	&158	&31	&108	&20	&238	&15\\
    9:30:17	&60	    &25	&24	    &7	&98	  &10\\
    9:30:18	&194	&19	&149	&5	&265  &9\\
    9:30:19	&101	&13	&25	    &9	&90	  &13\\
    9:30:20	&338	&29	&126	&10	&270  &12\\
    9:30:21	&45	    &8	&43	    &6	&42	  &5\\
    9:30:22	&59	    &10	&43	    &6	&52	  &4\\
    9:30:23	&32	    &11	&13	    &6	&41	  &8\\
    9:30:24	&632	&103&4	    &3	&80	  &7\\
    9:30:25	&85	    &36	&86	    &16	&64	  &6\\
    9:30:26	&34	    &10	&2	    &2	&58	  &7\\
    9:30:27	&33	    &7	&6	    &6	&16	  &2\\
    9:30:28	&73	    &17	&1	    &1	&14	  &2\\
    9:30:29	&55	    &9	&3	    &3	&33	  &8\\
    9:30:30	&69	    &14	&10	    &4	&38	  &7\\
    ... &... &... &...&... &... &...\\
  \hline
  \hline

\end{tabular}
}
\label{table: phenomenon of multiple transactions}
\end{table}
For example, Table \ref{table: phenomenon of multiple transactions} displays the numbers of transactions and different prices for Apple, Cisco Systems, and Microsoft Corporation in the first 30 seconds of trading on November 4, 2016.
In the case of Apple, there were 867 different transaction prices for a total of 5387 transactions in this time interval, and the number of transactions per second ranged from 32 to 632, with an average of about 180 transactions per second.
However, there was no information regarding the ordering of the different prices occurring at each second.
Such multiple transactions occur throughout the trading day, not only at the opening of the market.
\begin{figure}[!htb]
\centering
\includegraphics[width=12 cm,height =8.5cm]{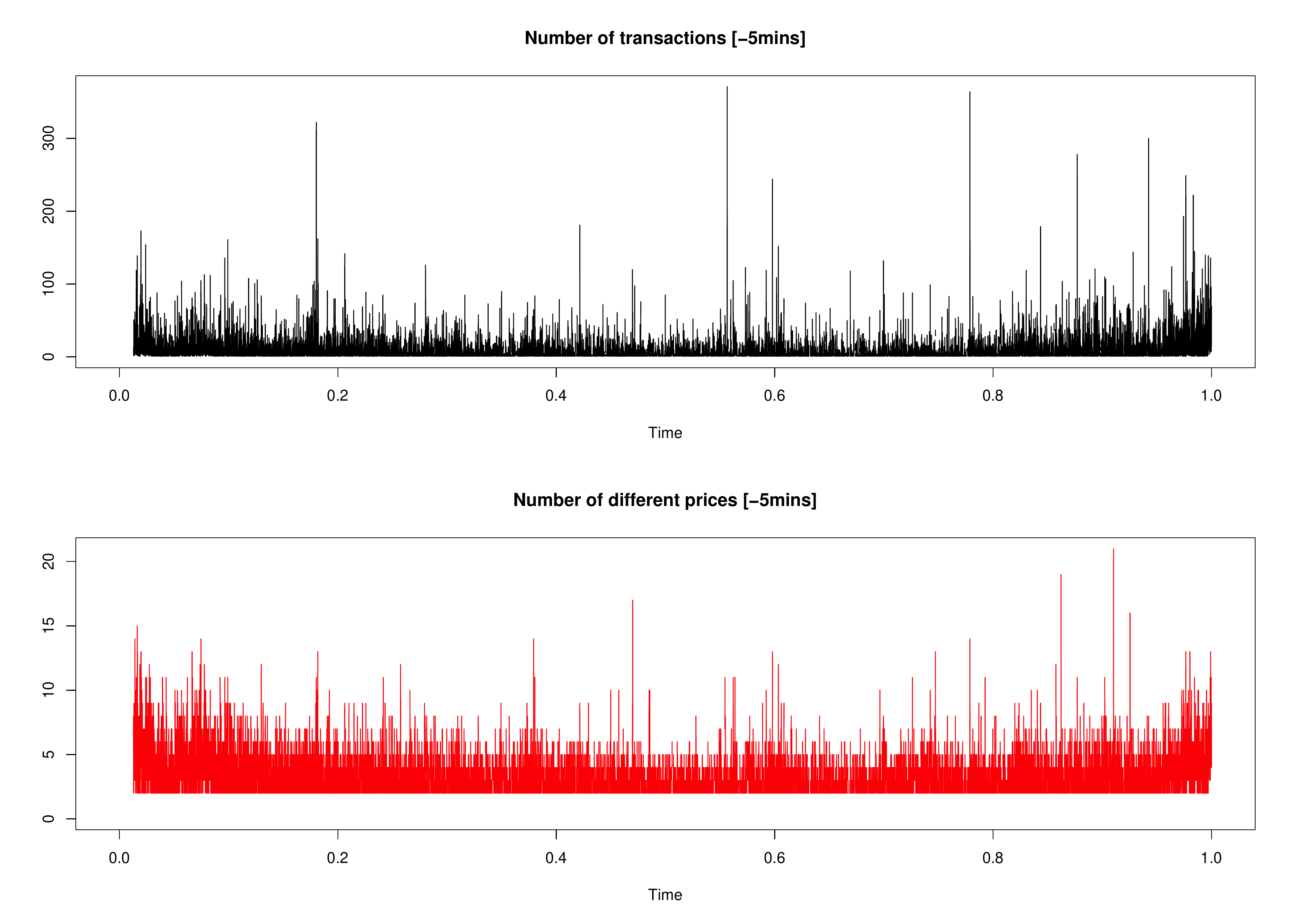}
\caption{ Number of transactions and number of different transaction prices of Apple stock during the first five minutes after opening of trading hours on November 4, 2016.
There were 161814 transactions for Apple stock during 23100 trading seconds, of which 11523 were multiple transactions (more than one transaction per second), approximately 49.88\% of the total transactions.}

\label{fig: transactions of apple in Nov4 2016}
\end{figure}
Figure \ref{fig: transactions of apple in Nov4 2016} shows in detail the numbers of transactions and different transaction prices of Apple stock during the first five minutes after opening of trading hours on November 4, 2016.
There were 161814 transactions for Apple stock, of which 11523 were multiple transactions (more than one transaction per second), representing approximately half of the total transactions.
In the presence of multiple transactions, the order of consecutive ticks is not observed and hence it is not possible to obtain the increment, as $\Delta \bX_i$ in RCV or TVA, at each time stamp.
Thus, it is natural to ask what effects multiple transactions have on the estimation of the ICV matrix, and whether an estimator of the ICV matrix can be obtained for use in the case of high-frequency multiple observations in a high-dimensional setting.
The effect of multiple transactions has been studied in the one-dimensional case;
see \cite{Jingetal2017}, \cite{Liu2017}, and \cite{Liuetal2018}.
However, few studies have focused on high-dimensional settings.
As the RCV matrix is not suitable for making a robust inference of the ICV matrix, the present work mainly focuses on investigating the limiting  behavior of a TVA matrix based on high-frequency data with multiple transactions.

\begin{figure}[!htb]
\centering
\includegraphics[width= 8cm,height =8cm]{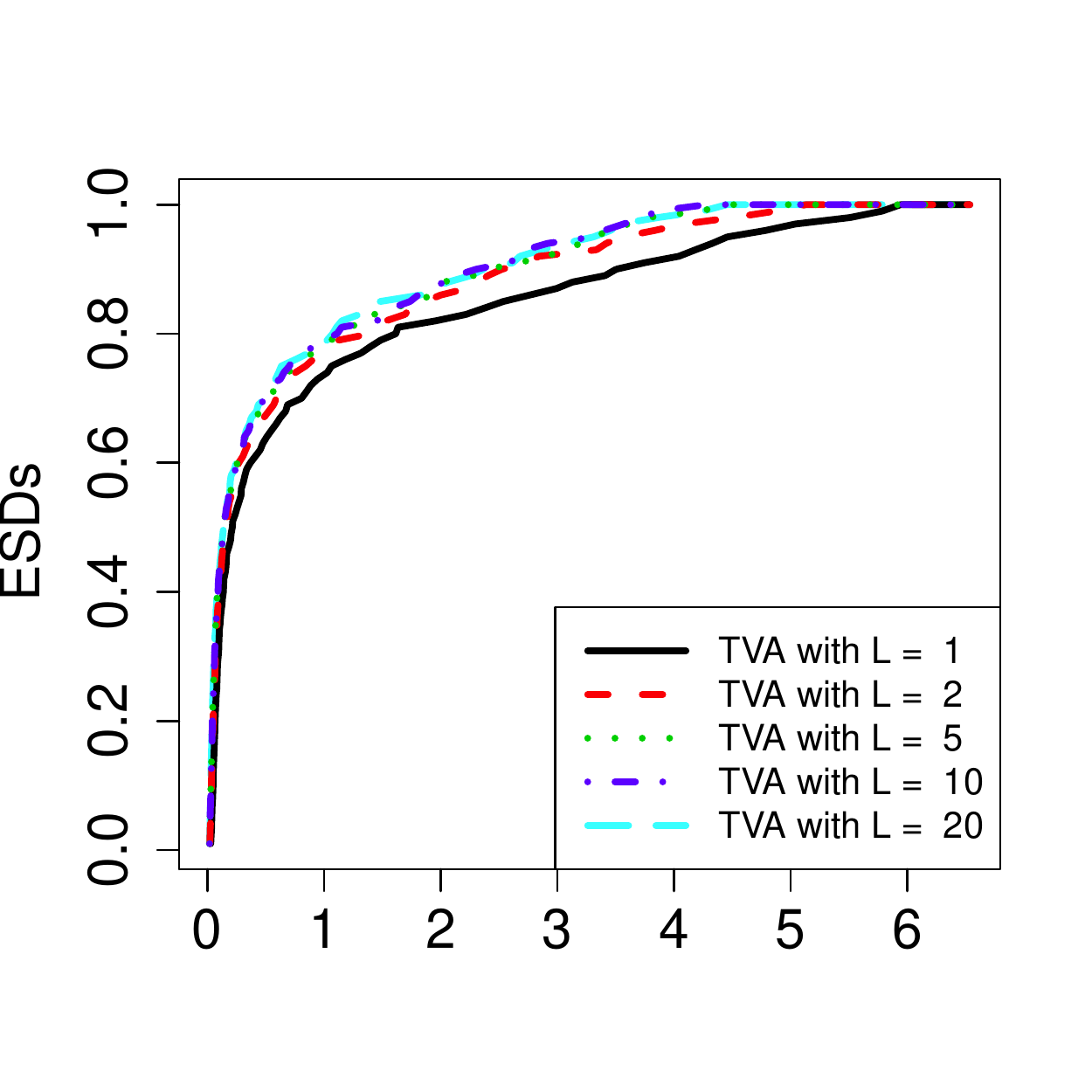}
\caption{
Plots of ESDs of TVA matrices with different values of $L$, for number of stocks $p=100$, observation frequency $n=390$, and multiple transaction frequency $L=1$, $2$, $5$, $10$, and $20$, respectively.
The log price processes $\{\bX_{s_j} \}$ are generated from the model $d\bX_s=\bLambda d\bW_s$ with $\bLambda=(0.5^{|i-j|})_{i,j=1,\cdots,p}$ and equally spaced time points $s_j=j/(nL)$ for $j=1,\cdots, nL$. The observations at each time $t_i=i/n$ are approximated by the averages of $L$ log prices, that is,
$ \bX_{t_i} = (1/L) \sum_{\ell=1}^L \bX_{s_{(i-1)L + \ell}}  $,
for $i=1,2,\cdots, n$. The TVA matrices were established based on observations $\{\bX_{t_i} \}$ for different values of $L$.
}
\label{fig: compare RCV TVA and their average}
\end{figure}

The effect of multiple transactions can be easily illustrated in the simplest case as follows.
We set the number of stocks to $p=100$, the
observation frequency to $n=390$ (one record per minute), and assume
equally spaced observation times ($t_i=i/n$) and an equal number of transactions at each time point, i.e., $L_i\equiv L$, $i=1,2,\cdots,n$.
Suppose that the log price processes $\{\bX_{s_j} \}$ are generated from the model $d\bX_s = \bLambda d\bW_s$, with $\bLambda=(0.5^{|k-\ell|})_{k,\ell=1,\cdots,p}$
and equally spaced time points $s_j=j/(nL)$, for $j=1,\cdots,nL$.
The observations $\bX_{t_i}$ are approximated by the averages of $L$ multiple transaction prices at each time point $t_i$, that is,
for $i=1,2,\cdots,n$,
$
\bX_{t_i} = (1/L) \sum_{\ell=1}^L \bX_{s_{(i-1)L + \ell}}.
$
From random matrix theory, for any $p$-dimensional symmetric matrix $\bA$, the empirical spectral distribution (ESD) is defined as
\begin{eqnarray}\label{eqn:ESD}
F^{\bA} (x) = \dfrac 1p \sum_{j=1}^p \ \indic (\lambda_j(\bA) \le x), ~~~~ \forall x\in \mathbb{R},
\end{eqnarray}
where $\indic(\cdot)$ is an indicator function and $\lambda_j(\bA)$ denotes the $j$th largest eigenvalue of matrix $\bA$.
The limit of ESD is referred to as the limiting spectral distribution (LSD), if it exists.
To illustrate the effects of multiple transactions, Figure \ref{fig: compare RCV TVA and their average}
presents the ESDs of TVA matrices with different values of $L$.
Clearly, the ESDs of TVA matrices depend  heavily on the value of $L$, the number of transactions,
and tend to be stable when $L$ is large.
This implies that the result of Theorem 2 in \cite{ZhengandLi2011} may no longer hold when multiple transactions occur.
Hence, it is necessary to develop new theoretical results for TVA matrices  and  propose an estimator of the ICV matrix based on high-frequency multiple observations.

In this work, we first study the limiting behavior of TVA matrices based on high-frequency multiple observations without considering microstructure noise.
We show that in the high-dimensional case, for a class $\mathcal{C}$ of diffusion processes, the LSD of an averaged TVA (ATVA) matrix depends not only on that of ICV but also on the numbers of multiple transactions at each recording time. Hence, we propose an adjusted ATVA (A-ATVA) estimator based on latent multiple high-frequency data.
Then, considering that
observed prices are always contaminated by market microstructure noise in practice, the limiting behavior of TVA matrices is further studied based on high-frequency noisy observations with multiple transactions.
We adopt the pre-averaging approach proposed in \cite{Jacodetal2009}
to deal with microstructure noise and
prove that for processes in class $\mathcal{C}$, the pre-averaging ATVA (PA-ATVA) matrix eliminates the effects of both  microstructure noise and multiple transactions, and its LSD depends solely on that of the ICV matrix.
To develop a consistent estimator of the ICV matrix, three types of nonlinear shrinkage (NS) estimators are considered  based on high-frequency noisy multiple observations.
Simulation studies are used to support our theoretical results and demonstrate the finite sample performance of the   proposed estimators.
Finally, high-frequency portfolio strategies are  evaluated under these estimators in real data analysis.

We give some notation that will be used throughout this article.
Let $\mathbb{C}$, $\mathbb{R}$, $\mathbb{Z}$, and $\mathbb{N}$ denote the sets of complex, real, integer, and natural numbers, respectively;  $\mathbb{C}^+$, a subset of $\mathbb{C}$, contains positive imaginary parts.
All vectors are column vectors, and we use $|\cdot|$ to denote the Euclidean norm for vectors.
The transpose of any matrix $\bA$ is denoted by $\bA^{\mT}$.
For any real matrix $\bA$, $\|\bA\| = \sqrt{\lambda_{\max}(\bA\bA^{\mathrm{T}})}$ denotes its spectral norm, where $\lambda_{\max}$ is the largest eigenvalue.
$||\bA||_F = \sqrt{\tr(\bA\bA^{\mT})}$ denotes the Frobenius norm of matrix $\bA$.
For any distribution function $F$, $m_F(\cdot)$ denotes its Stieltjes transform defined as
\[m_F(z) = \int \frac{1}{\lambda - z}dF(\lambda), \hspace{1cm} \text{for~} z\in \mathbb{C}^+ :=\{z \in \mathbb{C}: \Im(z) > 0\},\]
where  $\Im(z)$ is the imaginary part of $z$, and
$\lfloor\cdot\rfloor$ indicates rounding down to the nearest integer.

The rest of the paper is organized as follows.
The main theoretical results are presented in Section \ref{sec:results}, in which two asymptotic properties of ATVA matrices are established based on latent log price process and noisy observations with multiple transactions, respectively.
Further, three types of NS estimators are given in Section \ref{sec: nolinear shrinakge estimators}.
Simulation studies and real data analysis are shown in Section \ref{sec: simulation and Empirical Study}.
Conclusions are drawn in Section \ref{sec:conclusion}, and detailed proofs are provided in the Supplementary material.

\section{Main results}\label{sec:results}
The following model settings are assumed in the case of high-frequency data with multiple transactions.
In a given time interval $[0,1]$, one trading day, for any process $(\bV_t)$,
suppose that $(\bV_t)$ can be observed at time points $t_i = i/n$ and there are $L_i (L_i\geq 1)$ transactions during each time interval $(t_{i - 1}, t_i]$, for $i = 1,\cdots, n$.
Let $V_{T_{i-1}+j}$ be the $j$th observation at  transaction time $s_{T_{i-1}+j}$ during
time interval $(t_{i-1}, t_i]$ with $T_0\equiv 0$ and $T_i=\sum_{k=1}^i L_k$, for $j=1,2,\cdots,L_i$ and $i=1,2,\cdots, n$.
That is, $L_i$ transactions at time points $\{s_{T_{i-1}+j}, j=1,2,\cdots,L_i \}$ occurred during $(t_{i-1},t_i]$, where
$t_{i-1} < s_{T_{i-1}+1} < s_{T_{i-1}+2} <
\cdots < s_{T_{i-1}+L_i}=s_{T_i} \le t_i$.
The observations at each recording time point $\{t_i\}$                      are as follows:
\begin{description}
\item[] at time $t_0$: ~~ $\bV_{t_0} = \bV_0$;
\item[] at time $t_1$: ~~ $\bV_{1},\bV_{2},  \cdots, \bV_{T_1}$;
\item[] ~~$\vdots$
\item[] at time $t_n$: ~~ $\bV_{T_{n-1}+1}, \bV_{T_{n-1}+2}, \cdots \bV_{T_n}$.
\end{description}
Theoretically, these observations occur consecutively during time interval $(t_{i-1}, t_i]$. However, in practice, the observations at each recording time $t_i$ lose their order of arrival owing to the recording mechanism, as
is clear from Table \ref{table: phenomenon of multiple transactions}.
The $L_i$ observations occurring during time interval $(t_{i-1}, t_i]$ are only recorded at time $t_i$. Specifically,
\[\bV_0, \underbrace{\bV_{1}, \cdots, \bV_{T_1}}_{L_1 \ \text{observations at ~}t_1}, \cdots, \underbrace{\bV_{T_{i-1}+1}, \cdots, \bV_{T_i}}_{L_i \ \text{observations at ~}t_i}, \cdots, \underbrace{\bV_{T_{n-1}+1}, \cdots, \bV_{T_n}}_{L_n \ \text{observations at ~}t_n}.\]
One commonly used method to deal with multiple transactions involves taking averages of the multiple transaction prices at each time point $t_i$, that is,
\[
\ol{\bV}_i = \dfrac{1}{L_i}
\sum_{j=1}^{L_i} \bV_{T_{i-1}+j}, ~~~~~
i=1,\cdots,n,
\]
and using $\ol{\bV}_i$ to approximate the observation at time point $t_i$.
Then the increment becomes
\begin{align}
\Delta\wbV_i :&= \wbV_i - \wbV_{i-1}\notag \\
& = \frac{1}{L_i}\sum_{j=1}^{L_i} \bV_{T_{i-1} + j} - \frac{1}{L_{i-1}}\sum_{j=1}^{L_{i-1}} \bV_{T_{i-2}+ j}\notag\\
& = \frac{1}{L_i}\sum_{j=1}^{L_i} (\bV_{T_{i-1}+j} - \bV_{T_{i-2}}) - \frac{1}{L_{i-1}}\sum_{j=1}^{L_{i-1}} (\bV_{T_{i-2} + j} - \bV_{T_{i-2}} )\label{eqn:Delta_V}\\
& =  \sum_{j=1}^{L_i}a_{i,j}\Delta_{i ,j}\bV + \sum_{j=1}^{L_{i-1}}b_{i-1,j}\Delta_{i - 1,j}\bV\notag,
\end{align}
where $a_{i,j} = 1-\frac{j-1}{L_i}$, $b_{i,j} = \frac{j-1}{L_i}$, and
\[
\Delta_{i,j}\bV := \bV_{T_{i-1}+j} - \bV_{T_{i-1}+j-1}, ~\text{for}~~ j=1,\cdots, L_i, ~ i=1,\cdots,n,
\]
which is an asymmetric triangular form of $\Delta_{i,j}\bV$.
\subsection{LSD of non-overlapping averaged TVA (ATVA) matrix} \label{subsec:LSD_no_noise}
To investigate the effect of multiple transactions, we start by studying the non-overlapping averaged TVA (ATVA)  matrix, $\cA_N$, which is defined as

\begin{equation}
\cA_N :=
\dfrac{\sum_{i=1}^N |\Delta \ol{\bX}_{2i}|^2}{N} \cdot
\sum_{i=1}^N \dfrac{\Delta\ol{\bX}_{2i} \Delta\ol{\bX}_{2i}^T}{|\Delta\ol{\bX}_{2i}|^2}
= \dfrac{\sum_{i=1}^N |\Delta \ol{\bX}_{2i}|^2}{p}
\cdot \wt{\bSigma},
\end{equation}
where $N=\lfloor n/2 \rfloor$, $\wt{\bSigma}=\dfrac pN \sum_{i=1}^N \dfrac{\Delta\ol{\bX}_{2i} \Delta\ol{\bX}_{2i}^T}{|\Delta\ol{\bX}_{2i}|^2} $.
Using the notation in \eqref{eqn:Delta_V}, the log return based on the averaged log prices becomes
\[
\Delta\ol{\bX}_{2i} = \sum_{j=1}^{L_{2i}} a_{2i,j} \Delta_{2i,j} \bX
+ \sum_{j=1}^{L_{2i-1}}b_{2i-1,j} \Delta_{2i-1,j} \bX,
\]
where $\Delta_{i,j} \bX =\bX_{T_{i-1}+j} - \bX_{T_{i-1}+j-1}$.

Note that there are two main differences between the $\Delta \ol{\bX}_{2i}$ given above and that in (2.6) of \cite{XiaandZheng2018}, even though they are both described by the weighted mean of log returns.
First, for the $\Delta \ol{\bX}_{2i}$ above, the averaged number taken within each time interval $(t_{i-1},t_i]$ varies owing to the different values of the transaction time.
Thus, $\Delta \ol{\bX}_{2i}$ is written in an asymmetric triangular form.
However, for the $\Delta \ol{\bX}_{2i}$ in (2.6) of \cite{XiaandZheng2018}, the same average number was taken; thus, $\Delta \ol{\bX}_{2i}$ was expressed as a symmetric triangular form.
Second, in this work, to handle multiple transactions at each recording time, the average numbers taken in $\{\Delta \ol{\bX}_{2i},i=1,\cdots,n\}$ can be either finite or infinite.
In practice, the average numbers taken during each time stamp are usually considered as finite numbers.
However, an average number that went to infinity was taken in  (2.6) of \cite{XiaandZheng2018}, in order to deal with microstructure noise.

We now state our assumptions.
\begin{compactenum}\setcounter{enumi}{1}
\item[]
\begin{compactenum}
\item\label{asm:class_c} For all $p$, $(\bX_t)$ is a $p$-dimensional process in class $\mathcal{C}$ for some drift process $\bmu_t = (\mu_t^{(1)},\cdots,\mu_t^{(p)})^{\mathrm{T}}$ and covolatility process $(\bTheta_t) = (\gamma_t\bLambda)$;
\item\label{asm:bdd_mu} there exists a $C_0<\infty$ such that for all $p$ and all $l = 1,\cdots,p$, $|\mu_t^{(l)}|\leq C_0$ for all $t\in[0,1)$ almost surely;
\item\label{asm:bdd_gamma} there exists a $0\leq \delta_1<1/2$ and a sequence of index sets $\mathcal{I}_p$ satisfying $\mathcal{I}_p \subset \{1,\cdots, p\}$ and $\#\mathcal{I}_p = O(p^{\delta_1})$ such that $(\gamma_t)$ may depend on $\bW_t$ but only on $(W_t^{(l)}:l\in \mathcal{I}_p)$; moreover, there exists a $C_1<\infty$ such that for all $p$, $|\gamma_t|\in(1/C_1, C_1)$ for all $t\in [0,1)$, almost surely;
\item\label{asm:bdd_sigma} there exists $C_2<\infty$ such that for all $p$ and all $l$, the individual volatilities $\sigma_t^{(l)} = \sqrt{(\gamma_t)^2\sum_{k=1}^p(\Lambda_{lk})^2}\in (1/C_2, C_2)$ for all $t\in [0,1]$ almost surely; in addition, $(\gamma_t)$ converges uniformly to a nonzero process $(\gamma_t^*)$ that is piecewise continuous with finitely many jumps almost surely;
\item\label{asm:bdd_ICV} there exists $C_3<\infty$ and $0\leq \delta_2 <1/2$ such that for all $p$, $||\text{ICV}|| \leq C_3p^{\delta_2}$ almost surely;
\item\label{asm:bdd_delta} the $\delta_1$ in \eqref{asm:bdd_gamma} and $\delta_2$ in \eqref{asm:bdd_ICV} satisfy that $\delta_1 + \delta_2 <1/2$;
\item\label{asm:LSD_bH} almost surely, as $p\to\infty$, the ESD of $\bbSigma = \bLambda\bLambda^{\mathrm{T}}$ converges in distribution to a probability distribution $\breve{H}$;
\item\label{asm:dim_ratio} $N = \lfloor n/2 \rfloor$ with $\lim_{p\to \infty} p/N \to \rho \in (0, \infty)$;
\item\label{asm:bdd_duration} the transaction durations are equally spaced in each time interval $(t_{i-1}, t_i]$, that is, $\Delta s_{i,j}: = s_{T_{i-1}+j} - s_{T_{i-1}+j-1} = 1/(nL_i)$ for $j = 1,\cdots, L_i$ and $i = 1,\cdots, n$;
\item\label{asm:moment_L}  $\{L_{1},\cdots, L_n\}$ is a sequence of independent positive integers independent of $(\bX_t)$ with $E(\frac{1}{L_i}) = f_{1, t_i}$ and $E(\frac{1}{L_i^2}) = f_{2, t_i}$ for $i = 1, \cdots, n$, where $f_{1, t}$ and $f_{2, t}$ are c\`{a}dl\`{a}g functions satisfying $\lim\limits_{n\to\infty} f_{1, t} = f_{1, t}^*$, $\lim\limits_{n\to\infty} f_{2, t} = f_{2, t}^*$ for $t\in [0,1]$.

\end{compactenum}
\end{compactenum}

Assumptions \eqref{asm:class_c}--\eqref{asm:bdd_delta} are reasonable for the underlying log price process proposed in \cite{ZhengandLi2011}.
In such a case, no sparsity assumption on the ICV matrix is needed, and the dependence between the covolatility process and the Brownian motion in Assumption \eqref{asm:bdd_gamma} allows the leverage effect to be captured.
Assumptions \eqref{asm:LSD_bH}--\eqref{asm:dim_ratio} are standard for use in high-dimensional settings in random matrix theory.
The numbers of multiple transactions are considered in Assumptions \eqref{asm:bdd_duration}--\eqref{asm:moment_L} and are assumed to be equal and (weakly) stationary.
Assumption \eqref{asm:bdd_duration} is only for theoretical purposes, as it is not required when microstructure noise is involved.

\begin{thm}\label{thm:LSD_TVA}
Suppose that Assumptions \eqref{asm:class_c}--\eqref{asm:moment_L} hold.
Then, almost surely, the ESDs of ICV and $\cA_N$ converge to probability distributions $H$ and $F^{\cA}$, respectively, where
\begin{equation}
H(x)=\breve{H}(x/\theta), ~~~ \textrm{for all } x\geq 0 ~ \textrm{with} ~ \theta=\int_0^1 (\gamma_t^*)^2 dt,
\label{eq: determine of H}
\end{equation}
and $F^{\cA}$ and $H$ are related as follows:
\[
m_{\cA}(z) = \int_{\tau\in \mathbb{R}} \frac{1}{\tilde{\theta}_f/\theta\cdot\tau(1 - \rho(1 + zm_{\cA}(z))) - z}dH(\tau), ~~~~ \textrm{for} ~ z\in \mathbb{C}^+,
\]
where  $\tilde{\theta}_f = \int_0^1 (\dfrac 13 + \dfrac 16 f_{2,s}^*)(\gamma_s^*)^2ds$ and $m_{\cA}(z)$ denotes the Stieltjes transform of $F^{\cA}$.
\end{thm}

\begin{rmk}
Note that if $(\bX_t)$ belongs to class $\mathcal{C}$ then the ICV matrix can be written as $\ICV = \int_0^1 \gamma_t^2 dt\cdot \bbSigma$, where $\bbSigma = \bLambda\bLambda^{\mathrm{T}}$.
When $L_i \equiv 1$, Theorem 2 of \cite{ZhengandLi2011} shows that the LSDs of $\wt{\bSigma}$ and $\bbSigma$ are related through the Mar\v{c}enko--Pastur equation (see \eqref{eq: converfence of proposition A1}), and $2\sum_{i=1}^N|\Delta \bX_{2i}|^2/p $ is a consistent estimator of $\theta=\int_0^1 (\gamma_t^*)^2 dt$.
Throughout the proof of Theorem \ref{thm:LSD_TVA} in this paper, the LSDs of $\wt{\bSigma}$ and $\bbSigma$ are related through the Mar\v{c}enko--Pastur equation (see Proposition \ref{prop: convergence of tilde sigma} in the appendix); 
 however,  $2\sum_{i=1}^N|\Delta \bX_{2i}|^2/p $ is no longer consistent with $\theta$ in a general setting when $L_i \equiv 1$ is not satisfied. 
Hence, a non-overlapping A-ATVA matrix, $\wt{\cA}_N$, is proposed in Corollary \ref{cor:LSD_ATVA}, the LSD of which depends solely on that of the ICV matrix.
\end{rmk}

\begin{rmk}
Moreover, observe that the stock intraday volatility tends to be U-shaped, which implies that the numbers of transactions at opening and closing of the market tend to be much larger than the numbers in the middle of the day.
Thus, if Assumption \eqref{asm:moment_L} is further taken to be piecewise constant, Corollary \ref{cor:LSD_ATVA} (see below) holds.
\end{rmk}
If we further suppose that $(f_{2, s}^*)$ is piecewise constant, that is, there exist constant numbers  $\eta_1,\cdots,\eta_K$ for $K<\infty$ such that
\[
f_{2, s}^* = \sum_{i=1}^K \eta_i  \mathrm{I}\{s\in(a_{i-1}, a_i]\},
\]
where $(0, 1] = \sum_{i =1}^K (a_{i-1}, a_i]$ and $0 = a_{0} <a_1<\cdots<a_K = 1$, then we define the adjusted ATVA (A-ATVA) matrix as
\[
\wt{\cA}_N := \frac{1}{p}\sum_{i = 1}^K \left( \dfrac13 + \frac{1}{6(\ell_i - \ell_{i-1})}\sum_{j = \ell_{i-1}+1}^{\ell_i}\frac{1}{L_j^2}\right)^{-1} \sum_{m = \lfloor (\ell_{i-1} + 1)/2\rfloor + 1}^{\lfloor\ell_i/2\rfloor}|\Delta\wbX_{2m}|^2
\cdot \wt{\bSigma},
\]
where $\ell_i = \max\limits_{0\leq j\leq n}j/n\leq a_i$, provided that $\min\limits_{1\leq i \leq K}n(a_i - a_{i-1})\geq 3$.
A direct consequence of Theorem \ref{thm:LSD_TVA} is the following corollary, which demonstrates that the LSD of the A-ATVA matrix $\wt{\cA}_N $ depends solely on that of the ICV through the Mar\v{c}enko--Pastur equation.
\begin{cor}\label{cor:LSD_ATVA}
Under the assumptions of Theorem \ref{thm:LSD_TVA}, if we further suppose that $(f_{2, s}^*)$ is piecewise constant, then, almost surely, the ESD of $\wt{\cA}_N $ converges to a probability distribution $F^{\wt{\cA}}$, which is determined by $H$ through the Stieltjes transform via the Mar\v{c}enko--Pastur equation
\[
m_{\wt{\cA}}(z) = \int  \frac{1}{\tau(1 - \rho(1 + zm_{\wt{\cA}}(z))) - z}dH(\tau), ~~~~ \textrm{for} ~ z\in \mathbb{C}^+,
\]
where $m_{\wt{\cA}}(z)$ is the Stieltjes transform of $F^{\wt{\cA}}$.
\end{cor}


\subsection{Pre-averaging averaged TVA (PA-ATVA) matrix in the presence of microstructure noise}
\label{section: PA-TVA}
Besides multiple observations, microstructure noise presents another challenge to estimation of the ICV matrix.
It is commonly believed that, in practice, latent asset prices are always contaminated by market microstructure effects, so-called microstructure noise, which is induced by various frictions in the trading process such as asymmetric information among traders and bid-ask spread.
The high-frequency accumulation of microstructure noise seriously affects the influence of the ICV. 
Several methods have been developed in recent work to deal with microstructure noise, including the two/multi-scale approaches introduced by \cite{AitSahalietal2005}, \cite{MyklandandZhang2009}, \cite{AitSahalietal2010}, \cite{AitSahalietal2011}, \cite{Zhangetal2005}, and \cite{Zhang2006}; the realized kernel suggested by \cite{BNetal2008} and \cite{BNetal2011}; the quasi-maximum likelihood method studied by \cite{Xiu2010}; and the pre-averaging approach proposed by \cite{Jacodetal2009} and \cite{Li2013}.
In this section, we propose a PA-ATVA matrix to infer the ICV matrix, taking into consideration both multiple observations and microstructure noise.
Instead of the latent process $(\bX_t)$, the observed contaminated process $(\bY_t)$ is usually assumed to follow the additive model
\begin{equation}
\bY_{t} = \bX_{t} + \bm\varepsilon_{t}, \hspace{1cm} \text{for~} t\in [0,1],
\label{eq: model for microstructure noise}
\end{equation}
where $\bm\varepsilon_{t} = (\varepsilon_{t}^{(1)}, \cdots,\varepsilon_{t}^{(p)})^\mathrm{T}$ denotes the noise process.
In the presence of multiple observations, the noisy observations $(\bY_t)$ occur as follows:
\[\underbrace{\bY_{1}, \cdots, \bY_{T_1}}_{L_1\text{~observations at ~}t_1}, \cdots, \underbrace{\bY_{T_{i-1}+1}, \cdots, \bY_{T_i}}_{L_i\text{~observations at ~}t_i}, \cdots, \underbrace{\bY_{T_{n-1}+1}, \cdots, \bY_{T_n}}_{L_n\text{~observations at ~}t_n}.\]
For any process $\bV_t$, recall that the average of multiple observations at each recording time $t_i$ is denoted by
\[
\ol{\bV}_i = \dfrac{1}{L_i}
\sum_{j=1}^{L_i} \bV_{T_{i-1}+j}, ~~~~~
i=1,\cdots,n.
\]
We further introduce the following notation:
\[  \hspace{1cm}\twbV_i:= \frac{1}{h}\sum_{j = 1}^{h}\widebar{\bV}_{(i-1)h+j}  ~~~\text{and}~~~\Delta\twbV_{2i} = \twbV_{2i} - \twbV_{2i - 1}.\]
In this work, we adopt a pre-averaging method to deal with microstructure noise. Then, using the notation introduced above, the observed return based on the pre-averaged price becomes
\[
\Delta\twbY_{2i} = \Delta\twbX_{2i} + \Delta\twbe_{2i}.
\]
Next, we define an averaged version of the TVA matrix as follows.
Let $h = \lfloor\xi n^{\beta}\rfloor$ with $\xi \in (0,\infty)$ and $\beta \in (1/2,1)$.
Take $\hv = \lfloor n/(2h)\rfloor $.
The PA-ATVA matrix is then defined as
\begin{equation*}
\mbB_{\hv}:= 3\frac{\sum_{i=1}^\hv|\Delta\twbY_{2i}|^2}{\hv}\cdot \sum_{i=1}^\hv\frac{\Delta\twbY_{2i}(\Delta\twbY_{2i})^{\mathrm{T}}}{|\Delta\twbY_{2i}|^2} =  3\frac{\sum_{i=1}^\hv|\Delta\twbY_{2i}|^2}{p}\cdot\hbSigma,
\end{equation*}
where
\begin{equation}
\hbSigma := \frac{p}{\hv} \sum_{i=1}^\hv\frac{\Delta\twbY_{2i}(\Delta\twbY_{2i})^{\mathrm{T}}}{|\Delta\twbY_{2i}|^2}.
\label{eq: definition of hbSigma}
\end{equation}
One key observation is that the window length $h$ has a higher order than $\sqrt{n}$, which enables us to asymptotically eliminate the effect of microstructure noise.
To study the behavior of $\mbB_{\hv}$  based on noisy and multiple observations, we require some assumptions regarding the noise process.
Recall that the definition of $\rho$-mixing coefficients is as follows.
\begin{defn} \label{def: rho-misiing}
For a stationary time series $(U_k),k\in \mathbb{Z}$, let $\mathcal{F}_j^\ell$ be the $\sigma$-field generated by the random variables $(U_k: -\infty\leq j\leq k \leq \ell\leq\infty)$.
The $\rho$-mixing coefficients are defined as
\[\rho(r) = \sup\limits_{f\in\mathcal{L}^2(\mathcal{F}_{-\infty}^0) ,~~ g\in \mathcal{L}^2(\mathcal{F}_{r}^\infty)} |\begin{rm}Corr\end{rm}(f,g)|, \hspace{0.5cm} \text{for} \hspace{0.5cm} r\in\mathbb{N},\]
where, for any probability space $\Omega$, $\mathcal{L}^2(\Omega)$ refers to the space of square-integrable, $\Omega$-measurable random variables.
\end{defn}
\begin{compactenum}\setcounter{enumi}{2}
\item[]
\begin{compactenum}
    \item\label{asm:error} For all $j= 1,\cdots,p$, the noise $(\varepsilon_t^{(j)})$ is a stationary time series with mean 0 and has bounded $4\ell$th moments and $\rho$-mixing coefficients $\rho^{(j)}(r)$ satisfying $\max_{j=1,\cdots,p} \rho^{(j)}(r) = O(r^{-\ell})$ for some integer $\ell>2$;
    \item\label{asm:dim_windows} $h =\lfloor \xi n^{\beta} \rfloor$ for some $\xi\in(0,\infty)$ and $\beta \in ((3+\ell)/(2\ell +2), 1)$, and $\hv = \lfloor n/(2h)\rfloor$ satisfy that $\lim_{p\to\infty} p/\hv = c >0$, where $\ell$ is the integer in Assumption \eqref{asm:error};
    \item\label{asm:bdd_L}  there exists a constant $L^*<\infty$ such that  $\sup\limits_{1\leq i \leq n} L_i\leq L^*$, almost surely;
    \item\label{asm:interval}   $\max\limits_{1\leq i \leq n} n(s_{T_i} - s_{T_{i-1}}) \to 1$, almost surely, as $n\to \infty$.
\end{compactenum}
\end{compactenum}

As pointed out by \cite{XiaandZheng2018}, Assumption \eqref{asm:error} is a quite `mild' assumption that allows for not only dependence within the noise process, both cross-sectional and temporal, but also dependence between the noise and price processes.
We have the following convergence result for the PA-ATVA matrix $\mbB_M$.

\begin{thm}
Suppose that Assumptions \eqref{asm:class_c}-\eqref{asm:LSD_bH} and \eqref{asm:error}-\eqref{asm:interval} hold.
Then, as $p\to\infty$, the ESDs of ICV and $\mbB_M$ converge almost surely to probability distributions $H$ and $F^{\mbB}$, respectively, where $H$ satisfies \eqref{eq: determine of H} and $F^{\mbB}$ is determined by $H$ in that its Stieltjes transform, $m_{\mbB}(z)$, satisfies the Mar\v{c}enko--Pastur equation,
\begin{equation}\label{eq: m-p equation}
m_{\mbB}(z) = \int_{\tau\in \mathbb{R}} \frac{1}{\tau(1 - c(1 + zm_{\mbB}(z))) - z}dH(\tau),\hspace{0.2cm}\text{for}~z\in \mathbb{C^{+}}.
\end{equation}
\label{thm: theorem for the TVAPA}
\end{thm}
In addition to microstructure noise, another challenge of high-frequency data analysis is asynchronous trading.
In practice, different stocks have different numbers of transactions during one time stamp.
Take the example of Apple, Cisco Systems, and Microsoft in Table \ref{table: phenomenon of multiple transactions}; there were 491, 166, and 217 transactions in the first trading seconds on November 4, 2016, respectively.
Further, denote by $L_i^{(q)}$ the number of multiple transactions for stock $q$, $1\leq q \leq p$, during recording interval $ (t_{i-1}, t_i]$.
In financial markets, different stocks commonly have different values of $L_i^{(q)}$.
However, one big breakthrough of the current work is that Theorem \ref{thm: theorem for the TVAPA} still holds even when $L_i^{(q)}$ are different for different values of $q$.
For any process $(\bV_t)$, let $V_{i,j}^{(q)}$ denote the observation of the $j$th transaction for stock $q$ during time interval $(t_{i-1}, t_i]$.
The true transaction time of $V_{i,j}^{(q)}$ is denoted as $s_{T_{i-1}+j}^{(q)}$, for $j = 1,\cdots, L_i^{(q)}$, satisfying $t_{i-1} \leq s_{T_{i-1}}^{(q)}< s_{T_{i-1} + 1}^{(q)}<\cdots < s_{T_{i-1} + L_i^{(q)}}^{(q)} = s_{T_i}^{(q)} \leq t_i$.
With asynchronous trading, the average of multiple observations at each recording time $t_i$ is denoted as
\[
\ol{\bV}_i^* = \left(\sum_{j=1}^{L_i^{(1)}}\frac{1}{L_i^{(1)}}V_{i,j}^{(1)}, \cdots, \sum_{j=1}^{L_i^{(p)}}\frac{1}{L_i^{(p)}}V_{i,j}^{(p)} \right)^{\mT},
\] 
and the pre-averaging averaged observation becomes
\[
\twbV_i^* = \frac{1}{h}\sum_{j = 1}^{h}\widebar{\bV}_{(i-1)h+j}^*,
\]
for $i = 1,\cdots, n$.
Using the above notation, the PA-ATVA matrix is then rewritten as
\[
 \mbB_M^* = 3\frac{\sum_{i=1}^\hv|\Delta\twbY_{2i}^*|^2}{\hv}\cdot \sum_{i=1}^\hv\frac{\Delta\twbY_{2i}^*(\Delta\twbY_{2i}^*)^{\mathrm{T}}}{|\Delta\twbY_{2i}^*|^2}.
\]
\begin{thm}\label{thm: theorem for the TVAPA asynt}
Suppose that Assumptions \eqref{asm:class_c}--\eqref{asm:LSD_bH} and \eqref{asm:error}--\eqref{asm:bdd_L} hold.
Further, assume that $\max\limits_{1\leq i \leq n, 1\leq q\leq  p} nh(s_{T_i}^{(q)} - t_i) \to 0 $, almost surely.
Then, as $p\to\infty$, the ESDs of ICV and $\mbB_M^*$ converge almost surely to probability distributions $H$ and $F^{\mbB^*}$, respectively, where $H$ satisfies \eqref{eq: determine of H} and $F^{\mbB^*}$ is determined by $H$ in that its Stieltjes transform $m_{\mbB^*}(z)$ satisfies the Mar\v{c}enko--Pastur equation,
\begin{equation}
m_{\mbB^*}(z) = \int_{\tau\in \mathbb{R}} \frac{1}{\tau(1 - c(1 + zm_{\mbB^*}(z))) - z}dH(\tau),\hspace{0.2cm}\text{for}~z\in \mathbb{C^{+}}.
\end{equation}
\end{thm}

\section{Nonlinear shrinkage estimators of ICV}
\label{sec: nolinear shrinakge estimators}
The limiting spectral behavior of matrix $\mbB_M$ established in Theorem \ref{thm: theorem for the TVAPA} implies that $\mbB_M$ is still a bad estimator of ICV.
Hence, our next problem was how to build a good estimator of ICV based on high-frequency multiple noisy observations.
\cite{LedoitandWolf2012} propose a class of rotation-equivariant estimators that retain the eigenvectors of the sample covariance matrix but nonlinearly shrink its eigenvalues by pushing up the small ones and pulling down the large ones without any particular structure assumption regarding the population covariance matrix.
This proposal makes perfect sense, as \cite{BaiandYin1993} and \cite{Yinetal1998} have provided solid evidence that the extreme eigenvalues of sample covariance matrices are more extreme than the population ones.
In this section, motivated by the idea of an NS strategy, we adjust the method proposed in \cite{LedoitandWolf2012} and further propose three types of NS estimators of the ICV matrix based on self-normalized noisy observations $\Delta \twbY_{2i}/|\Delta \twbY_{2i}|$.
A similar estimation strategy is applied for asynchronous trading data $\{Y_{i,j}^{(q)}\}$.

\subsection{(Original) nonlinear shrinkage (NS) estimator}\label{subsection: NS}
For any process in class $\cC$, note that $\ICV = \int_0^1 \gamma_t^2 dt \cdot \bbSigma$, in which $\bbSigma = \bLambda \bLambda^{\mT}$ is time invariant.
It was shown in the proof of Theorem \ref{thm: theorem for the TVAPA} that $\hat\theta_p = 3\sum_{i = 1}^M |\Delta \twbY_{2i}|^2/p$ is a consistent estimator of $\theta = \lim_{p\to\infty}\int_0^1 \gamma_t^2dt$, and the LSDs of $\wt{\bXi} = p/M \cdot \sum_{i = 1}^M\Delta \twbY_{2i} (\Delta \twbY_{2i})^{\mT}/|\Delta \twbY_{2i}|^2 $ and $\bbSigma$ are related through the Mar\v{c}enko--Pastur equation \eqref{eq: m-p equation}.
Suppose that we observe high-frequency multiple noisy observations $\{\bY_{t_i}\}$ during time period $[T - \tau, T]$ at current time $T$, for $\tau\geq 1$.
We construct matrix $\hbSigma$ based on observations $\{\bY_{t_i}\}$ over time period $[T - \tau ,T]$ and estimator $\hat\theta_p$ based on observations $\{\bY_{t_i}\}$ during time interval $[T - 1, T]$.
We denote the spectral decomposition of matrix $\hbSigma$ by
\[
\hbSigma = \bP_{T - \tau, T}\wt{\bLambda}_{T - \tau, T} \bP_{T - \tau, T}^{\mT},
\]
where the eigenmatrix $\bP_{T - \tau, T} = (\bp_{T - \tau, T}[1], \cdots, \bp_{T - \tau, T}[p])$, whose $i$th column is the eigenvector of $\hbSigma $, and $\wt{\bLambda}_{T - \tau, T} = \diag \{ \wt{\lambda}_1 ,\cdots,  \wt{\lambda}_p\}$ is a diagonal matrix with non-increasing eigenvalues $\wt{\lambda}_1\geq \cdots \geq \wt{\lambda}_p$.
To construct an NS estimator of $\bbSigma$, we retain the eigenmatrix $\bP_{T-\tau, T}$ and consider a class of estimators of the form:
\[
\bP_{T-\tau, T} {\bf D} \bP_{T-\tau, T}^{\mT},
\]
where $D = \diag( d_1, \cdots,  d_p)$ is diagonal.
In order to measure the distance between two matrices, we adopt the Frobenius norm.
The objective is to find the diagonal matrix ${\bf D}$ that solves the optimization problem
\[
\min\limits_{{\bf D} ~~\text{diag}}\| \bP_{T - \tau, T} {\bf D} \bP_{T - \tau, T}^{\mT} - \bbSigma\|_F.
\]
By elementary calculation, the solution is
\[D^{or} = \diag(d_1^{or},\cdots, d_p^{or}),
\]
where
\begin{equation}\label{eq: oracle di}
d_i^{or} = (\bp_{T - \tau, T}[i])^{\mT} \bbSigma \bp_{T - \tau, T}[i].
\end{equation}
Therefore, we propose the following nonlinear shrinkage (NS) estimator for ICV:
\[
\wh{\bSigma}_{NS}: = \hat{\theta}_p \cdot \wh{\bXi}_{NS},
\]
where $\wh\bXi_{NS} = \bP_{T - \tau, T} \diag(\wh{d_1^{or}},\cdots, \wh{d_p^{or}}) \bP_{T-\tau ,T}^{\mT}$,
and $\wh{d_i^{or}}$ is an estimator of $d_i^{or}$, solved by the algorithm proposed by \cite{LedoitandWolf2017}.

\subsection{Averaged nonlinear shrinkage (ANS) estimator}
\label{subsection: ANS}
Owing to the difficulty and computational complexity of estimation of $d_i^{or}$ in \eqref{eq: oracle di}, \cite{Lam2016} obtained a nonparametric eigenvalue-regularized covariance estimator by splitting the data into two parts and regularizing the eigenvalues using two independent sets of split data.
The theoretical properties of the regularized eigenvalues in \cite{Lam2016} have been established, showing that they are asymptotically the same as those $d_i^{or}$ in \eqref{eq: oracle di}.
Further, Lam improved the estimator by taking an average of the permutation data.
The excellent finite sample performance of this estimator has been demonstrated in various settings.
In this section,  we adopt the data-splitting idea of \cite{Lam2016} and propose an ANS estimator for ICV matrices in the high-frequency setting.
Suppose that we have tick-by-tick high-frequency observations $\{\bY_{l_i}\}$ over time period $[T - \tau, T]$ at current time $T$.
First, the estimator $\hat{\theta}_p = 3\sum_{i=1}^M |\Delta \twbY_{2i}|^2/p$ can be established based on high-frequency data during time interval $[T - 1, T]$.
Suppose that there are in total $M_{\tau} = \tau \cdot M$ averaged noisy increments $\Delta \twbY_{2j}, j = 1, \cdots, M_{\tau}$, over time period $[T - \tau, T]$, in which the overnight increments are excluded and there are $M$ increments each day.
Next, we permute the increments $B$ times.
At the $k$th permutation, we split the increments $\Delta \twbY^{(k)}$ into two parts, say $\Delta \twbY^{(k)} = (\Delta \twbY^{(k)}_1, \Delta \twbY^{(k)}_2)$, for $k = 1,\cdots, B$, with $\Delta \twbY^{(k)}_\ell$ having size $p\times M_{\tau_{\ell}}$ $(\ell = 1, 2)$ and $M_{\tau_1} + M_{\tau_2} = M_{\tau}$.
We define
\[
\hbSigma_{\ell}^{(k)} := \frac{p}{\hv_{\tau_{\ell}}}\sum_{i \in J_{\ell, k}}\frac{\Delta\twbY_{2i} (\Delta\twbY_{2i})^{\mathrm{T}}}{|\Delta\twbY_{2i}|^2},
\]
where $J_{\ell, k} = \{i: \Delta \twbY_{2i} \in \Delta\twbY_{\ell}^{(k)}\}$ for $\ell = 1, 2$ and $k = 1, \cdots, B$.
We denote the spectral decomposition of $\wt{\bXi}_1^{(k)}$ as
\[
\wt{\bXi}_1^{(k)} = \bP_1^{(k)}\wt{\bLambda}_1^{(k)}(\bP_1^{(k)})^{\mT},
\]
for $k = 1,\cdots, B$.
Then the averaged nonlinear shrinkage (ANS) estimator of ICV is defined as
\begin{equation*}
\wh{\bSigma}_{\ANS}:= \hat{\theta}_p\cdot\frac{1}{B} \sum_{k=1}^B \hbSigma^{(k)}_{\ANS}, ~~\hbSigma^{(k)}_{\ANS} = \bP_1^{(k)} \diag((\bP_1^{(k)})^{\mT} \wt{\bXi}_2^{(k)}\bP_1^{(k)}) (\bP_1^{(k)})^{\mT},
\end{equation*}
where $\diag(\bA)$ sets all non-diagonal elements of matrix $\bA$ to zero.
To this end, $M_{\tau_1}$ is chosen using the following criterion,
\begin{equation}\arg\min\limits_{\hv_{\tau_1}\in \mathcal{G}} \bigg| \bigg|  \frac{1}{B} \sum_{k = 1}^B (\hbSigma^{(k)}_{\ANS} -  \hbSigma_2^{(k)}) \bigg| \bigg|_F^2,
\label{eq: nonlinear shrnkage criterion}\end{equation}
where
\begin{equation}
\begin{split}
\mathcal{G} = \bigg[2\hv_{\tau}^{1/2}, 0.2\hv_{\tau}, 0.4\hv_{\tau},
0.6\hv_{\tau}, 0.8\hv_{\tau},
 \hv_{\tau} - 2.5\hv_{\tau}^{1/2}, \hv_{\tau} - 1.5\hv_{\tau}^{1/2}\bigg]
\end{split}
\label{eq: critia for selection of n1}
\end{equation}
is the candidate set for $\hv_{\tau_1}$, recommended by \cite{Lam2016} and \cite{Lametal2017}.

\subsection{Mixed nonlinear shrinkage (MNS) estimator}
\label{subsection: MNS}
Note that the effective sample size drops dramatically from $n$ to $M = O(n^{1-\beta})$, for $\beta \in (1/2,1)$, when the matrix $\wt\bXi$ in \eqref{eq: definition of hbSigma} is constructed.
In fact, only diagonal entries of matrix $(\bP_1^{(k)})^{\mT} \wt{\bXi}_2^{(k)}\bP_1^{(k)}$ are needed in the estimator $\wh{\bSigma}_{\ANS}$.
Hence, motivated by the idea of \cite{Liuetal2016}, a third type of NS estimator $\wh{\bSigma}_{\MNS}$ is considered as follows.
Suppose that we observe high-frequency data $\{\bY_t\}$ with multiple transactions during time period $[T - \tau, T]$.
First, the matrix $\wt{\bXi}_{T-\tau, T-1}$ in \eqref{eq: definition of hbSigma} is constructed based on high-frequency data during time period $[T - \tau, T - 1]$, and its spectral decomposition is denoted as
\[
\wt{\bXi}_{T-\tau, T - 1} = \bP_{T-\tau, T - 1}\wt{\bLambda}_{T-\tau, T - 1}(\bP_{T-\tau, T - 1})^{\mT},
\]
where
\[
\bP_{T-\tau, T - 1} = (\bp_{T-\tau, T - 1}[1], \cdots, \bp_{T-\tau, T - 1}[p]).
\]
Second, we set the eigenvectors $\bp_{T-\tau, T - 1}[k]$, $k = 1,\cdots, p$, to be given beforehand, and transform the averaged high-frequency data $\{\ol{\bY}_i\}$ during time period $[T - 1, T]$.
That is, consider $\ol{\bY}_{k, i}^{**} = (\bp_{T-\tau, T - 1}[k])^{\mT}\ol{\bY}_i$ as the transformed averaged observations at recording time $t_i \in [T - 1,T]$, for $k = 1,\cdots, p, i = 0,\cdots, n$.
Third, for each $k = 1,\cdots, p$, compute the averaging-pre-averaging estimator 
\[
\wh{d_i^{\APA}} = \frac{12}{k_n }\sum_{i = 0}^{n - k_n + 1}(\Delta_{i,k_n}^k\ol{\bY}_{k, i}^{**} )^2 - \frac{6}{\vartheta^2 n} \sum_{i = 0}^{n-1} (\ol{\bY}_{k, i + 1}^{**} - \ol{\bY}_{k, i}^{**})^2,
\]
where $k_n = \lfloor \vartheta n^{1/2}\rfloor,  \vartheta>0$ and $\Delta_{i,k_n}^k\ol{\bY}_{k, i}^{**} = \sum_{j = 1}^{k_n -1} g(j/k_n)( \ol{\bY}_{k, i + j + 1}^{**} - \ol{\bY}_{k, i + j}^{**})$ with $g(x) = x\wedge (1-x)$, which is a consistent estimator of the integrated volatility of $(\ol{\bY}_{k, i}^{**})$ over $[T - 1, T]$ by given $\bU_{T - \tau, T - 1}$, see \cite{Jacodetal2009} and \cite{Liuetal2018}.
Finally, our mixed nonlinear shrinkage (MNS) estimator of ICV is defined as
\[
\wh{\bSigma}_{\MNS} = \bP_{T -\tau, T -1} \diag (\wh{d_1^{\APA}}, \cdots, \wh{d_p^{\APA}})\bP_{T -\tau, T -1}^{\mT}.
\]

\section{Simulation and Empirical Study}
\label{sec: simulation and Empirical Study}
In this section, we demonstrate the finite sample performances of Theorem \ref{thm:LSD_TVA} and Theorem \ref{thm: theorem for the TVAPA} by showing that the proposed estimators $\mbA_N$ and $\mbB_M$ have almost the same empirical spectral distributions as their corresponding sample covariance matrices based on independent and identically distributed (i.i.d.) samples drawn from the ICV matrix.
More impressively, the asynchronism of different stocks with different trading numbers at one time stamp does not disable the recommended estimator $\mbB_M$.
Also, we show that the three types of NS estimators perform very well even with spiked eigenvalues where the largest few eigenvalues differ from the rest.

We adopt scenarios for the diffusion process $(\bX_t)$ from \cite{XiaandZheng2018} and \cite{Lametal2017}.
We take the following U-shaped stochastic process $(\gamma_t)$ as
\begin{equation}
d\gamma_t = -\rho(\gamma_t - \mu_t)dt + \sigma d\widetilde{W}_t, \hspace{0.5cm}\text{for} \hspace{0.5cm} t\in[0,1],
\label{eq: the definition of gamma in simulation}
\end{equation}
where $\rho = 10, \sigma = 0.05, \mu_t = 2\sqrt{0.0009 + 0.0008\cos(2\pi t)}$, and the process $\widetilde{W}_t = \sum_{i=1}^p W_t^{(i)}/\sqrt{p} $ with $W_t^{(i)}$ being the $i$th component of the Brownian motion $(\bW_t)$ that derives the price process. 
\eqref{asm:bdd_gamma} is violated since $(\gamma_t)$ depends on all components of the Brownian motion.
However, our estimates still work, as demonstrated by the simulation studies.
We assume that $\bLambda = (0.5^{|i-j|})_{i ,j = 1,\cdots,p}$ and further rescale it to satisfy the condition $\tr(\bLambda \bLambda^{\mT}) = p$ when spiked eigenvalues are considered.
The latent log price process $(\bX_t)$ follows
\begin{equation}\label{eq: simulation model}
d\bX_t = \gamma_t\bLambda d\bW_t.
\end{equation}

\subsection{Simulation results for A-ATVA matrix}
\begin{figure}[!htb]
\centering
\includegraphics[width= 6cm, height = 6cm]{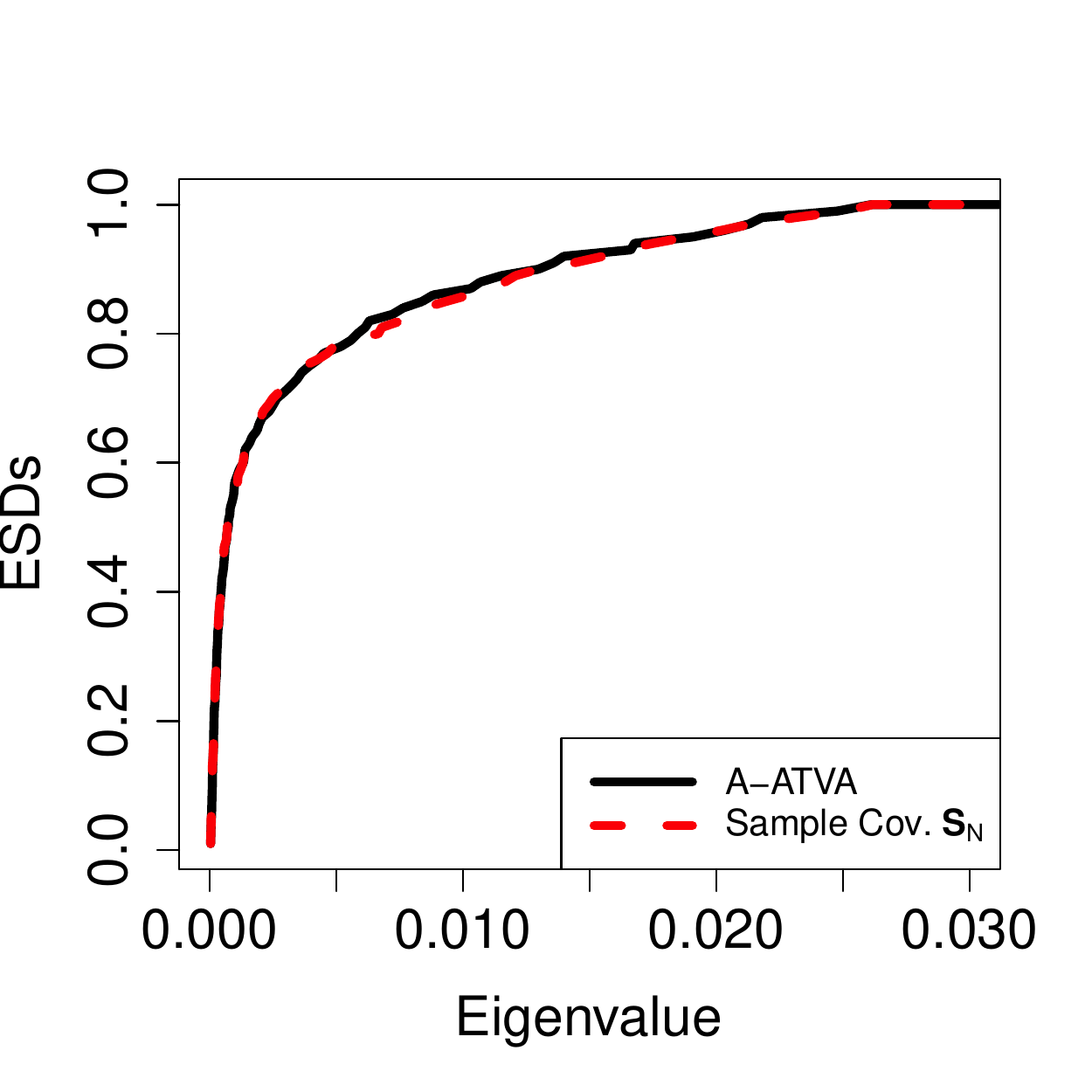}
\caption{ESDs of A-ATVA matrix $\mbA_{N}$  and sample covariance matrix $\bS_N$ for dimension $p = 100$ and observation frequency $n = 390$.
The number of multiple transactions $L_i -1$ for $i =1,\cdots, n$, is generated independently from a Poisson distribution with parameter 5.
}
\label{fig: application of theorem 1}
\end{figure}

\begin{figure}[!htb]
\centering
\includegraphics[width= 6cm, height = 6cm]{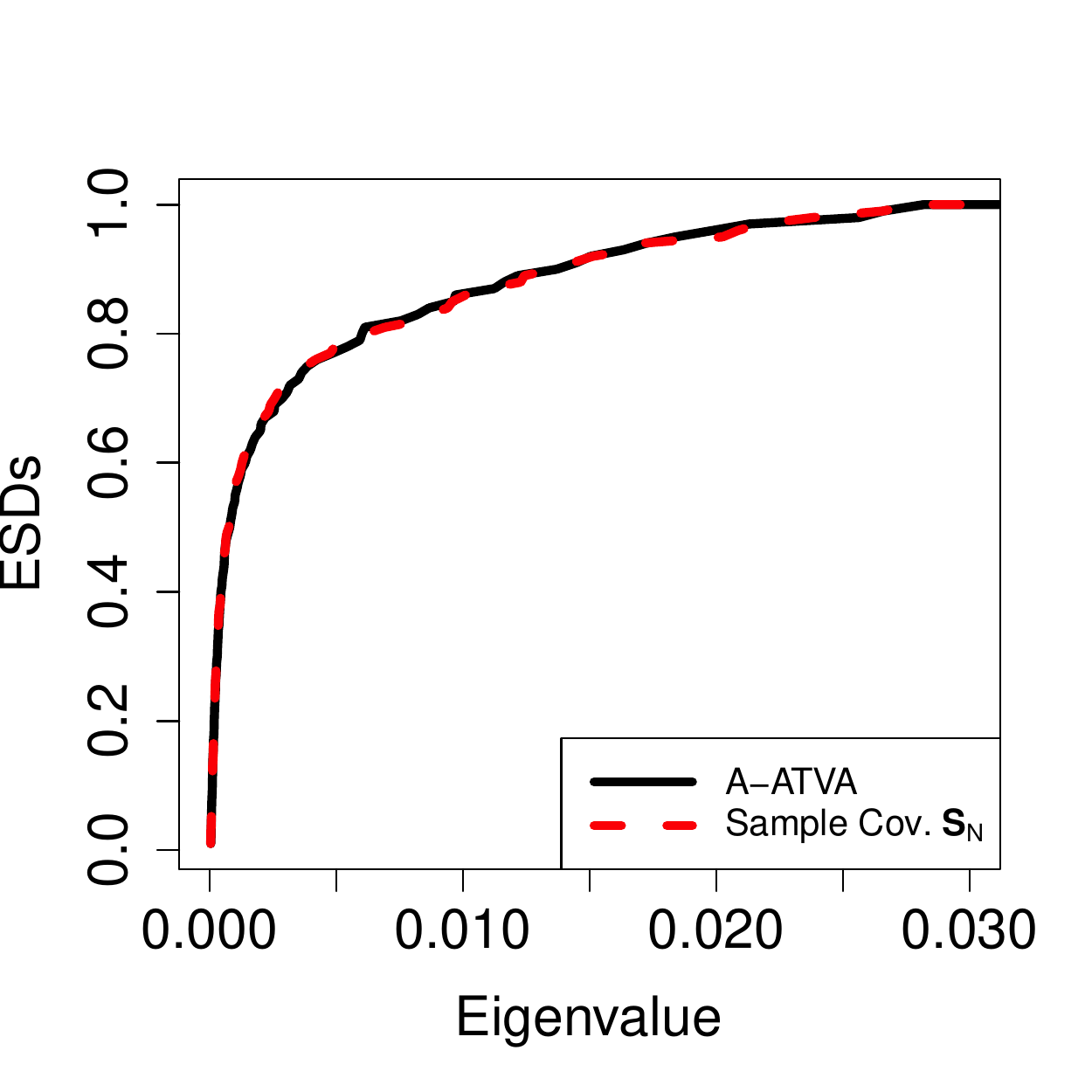}
\caption{ESDs of A-ATVA matrix $\mbA_{N}$  and sample covariance matrix $\bS_N$ for dimension $p = 100$ and observation frequency $n = 390$.
The numbers of multiple transactions $L_i -1$  are generated independently from a Poisson distribution with parameter 20 for the first and last $1/6$ of the data and with parameter 5 for the remaining $2/3$ of the data.
}
\label{fig: application of theorem 1 different distribution}
\end{figure}

\begin{figure}[!htb]
\centering
\includegraphics[width= 6cm, height = 6cm]{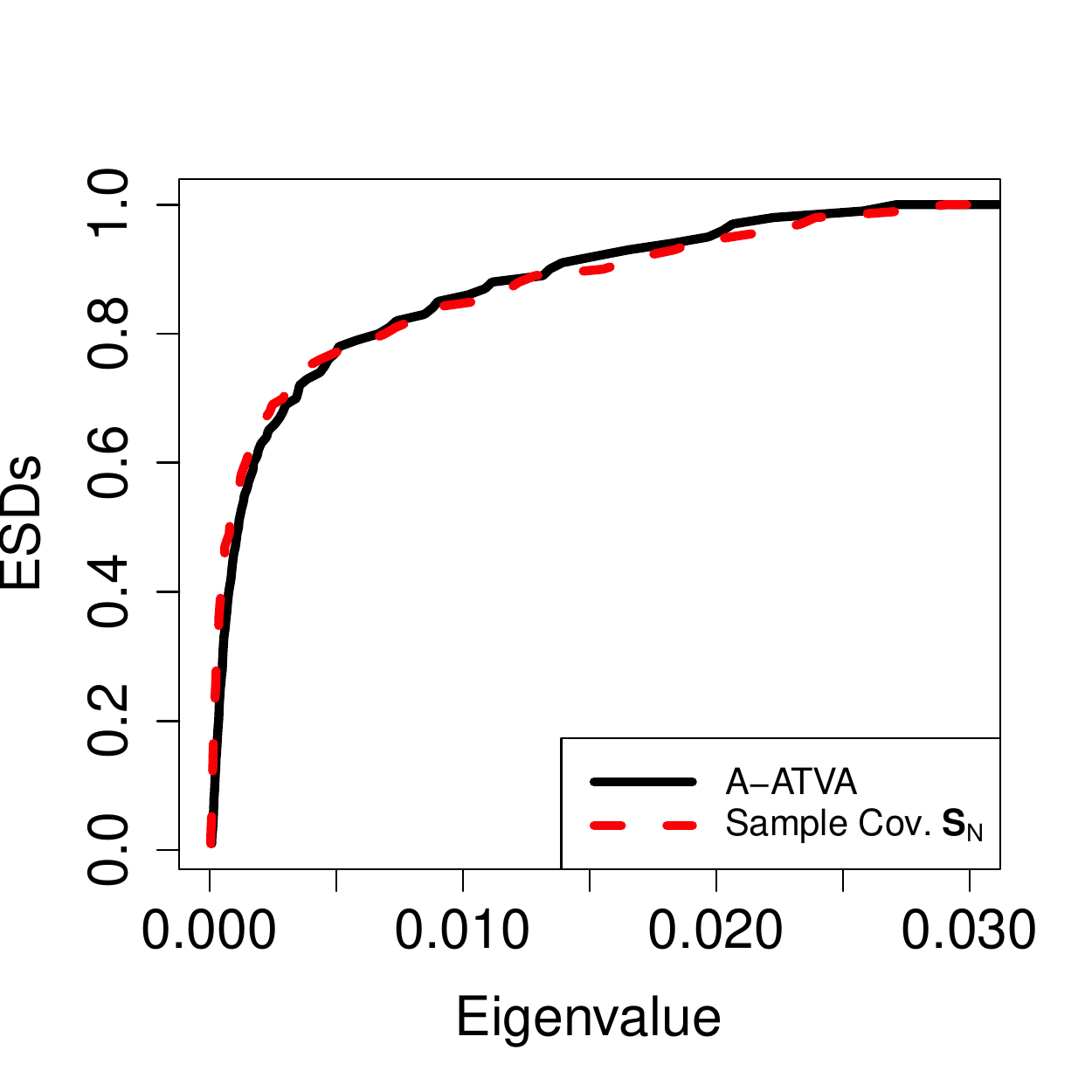}
\caption{ESDs of A-ATVA matrix $\mbA_{N}$  and sample covariance matrix $\bS_N$ for dimension $p = 100$ and observation frequency $n = 390$.
Let $L_i^{(q)}$ be the number of multiple transactions for stock $q$ within recording time stamp $(t_{i-1}, t_i]$, for $q = 1, \cdots, p, i = 1,\cdots, n$.
Then each  $L_i^{(q)}$  is generated independently from discrete uniform distribution $U[1,5]$.
}
\label{fig: application of theorem 1 ill}
\end{figure}

In this subsection, we compare the ESD of the A-ATVA matrix $\mbA_N$ and that of the sample covariance matrix $\bS_N$ with i.i.d. samples drawn from the ICV matrix, that is, $\bS_{N} = 1/N \sum_{i=1}^N (\ICV)^{1/2} \bZ_i \bZ_i^{\mT} (\ICV)^{1/2}$, where $\bZ_i \overset{i.i.d.}{\sim} N(0, \bI)$ and $N = \lfloor n /2\rfloor$.
It is well known that the LSDs of a sample covariance matrix $\bS_N$ and the underlying population covariance matrix ICV are related through the Mar\v{c}enko--Pastur equation from Theorem 1.1 of \cite{Silverstein1995}, which coincides with the result of Corollary \ref{cor:LSD_ATVA}.
We set $p = 100$ and $n = 390$, which represents the case where 100 stocks are recorded per minute within one trading day.
Figure \ref{fig: application of theorem 1} illustrates the simulation result when the number of multiple transactions $L_i - 1, i =1,\cdots, n$, are i.i.d. generated from a Poisson distribution with parameter 5, Poisson (5).
Intraday trading volatility commonly believed to be greater at the opening and closing of the market.
Figure \ref{fig: application of theorem 1 different distribution} shows the ESDs of the A-ATVA matrix $\mbA_N$ and sample covariance matrix $\bS_N$ when $L_i - 1$ are independently generated from Poisson (20) for the first and last $1/6$ of the data, and from Poisson (5) for the remaining $2/3$ of the data.
As shown in Figures \ref{fig: application of theorem 1} and \ref{fig: application of theorem 1 different distribution}, the ESDs of $\mbA_N$ and $\bS_N$ match quite well in both cases.
However, observe that different stocks usually have different numbers of transactions for each recording time stamp; see Table \ref{table: phenomenon of multiple transactions}.
The result of Corollary \ref{cor:LSD_ATVA} is further examined when $L_i^{(q)}$ (the number of transactions for stock $q$ during time interval $(t_{i-1}, t_i]$) varies for different $q = 1, \cdots, p$.
For each $i = 1,\cdots, n$ and $q = 1,\cdots, p$, we generate $L_i^{(q)}$ independently from a discrete uniform distribution, $U[1,5]$, and set the underlying trading time to be equally spaced during each recording time.
The resulting ESDs for the A-ATVA matrix $\mbA_N$ and sample covariance matrix $\bS_N$ are shown in Figure \ref{fig: application of theorem 1 ill}, and the maximum distances between two ESDs are reported in Table \ref{table: TVA fail discussion}, with $p/n = 10/39$ fixed but $n$ increasing.
\begin{table}[!htb]
\small
\centering
\setlength{\abovecaptionskip}{0pt}
\setlength{\belowcaptionskip}{10pt}
\caption{The averages of 1000 replications of maximum distance between the ESDs of A-ATVA matrix $\mbA_N$ and sample covariance matrix $\bS_N$ are reported with different values of $n = 390, 780, 1170$ but fixed p/n = 10/39.}
\setlength{\tabcolsep}{3 mm}{
\begin{tabular}{cccc}
  \hline
  \hline
   n  & 390 & 780 & 1170\\
   \hline
  distance & 0.13761 &0.13374 &0.13212\\
  \hline
  \hline
\end{tabular}
}
\label{table: TVA fail discussion}
\end{table}
The maximum distances between the two ESDs, $F^{\mbA_N}$ and $F^{\bS_N}$, can be estimated using the following formula
\begin{equation}\label{eq: estimation of maximum distance}
\max\limits_{x_1,\cdots,x_{2p}}|F^{\mbA_N}(x_i) -F^{\bS_N}(x_i) |,
\end{equation}
where $x_1,\cdots, x_{2p}$ are chosen to be equally spaced in the interval $[\min(\lambda(\mbA_N), \lambda(\bS_N) )$, $\max(\lambda(\mbA_N), \lambda(\bS_N) )]$, where $\lambda(\bA)$ denotes the spectral of any symmetric matrix $\bA$.
As shown in Figure \ref{fig: application of theorem 1 ill} and Table \ref{table: TVA fail discussion}, the distances between the ESDs of $\mbA_N$ and $\bS_N$ cannot be negligible when $L_i^{(q)}$ are different.
Thus, the limiting spectral of an A-ATVA matrix based on different values of $L_i^{(q)}$ remains an open problem. 

\subsection{Simulation results for PA-ATVA matrix}\label{subsec: PA-ATVA}
We now investigate the finite sample performance of the PA-ATVA matrix $\mbB_M$ in the presence of microstructure noise.
It is reasonable to conjecture that the ESD of the PA-ATVA matrix  $\mbB_M$ would have similar behavior to that of the sample covariance matrix $\bS_M = 1/M \sum_{i=1}^M (\ICV)^{1/2} \bZ_i \bZ_i^{\mT} (\ICV)^{1/2}$, where $\bZ_i \overset{i.i.d.}{\sim} N(0, \bI)$, as both the LSDs of
$\mbB_M$ and $\bS_M$ are related through the Mar\v{c}enko--Pastur equation in Theorem \ref{thm: theorem for the TVAPA} and Theorem 1.1 of \cite{Silverstein1995}.
Hence, the ESDs of the PA-ATVA matrix $\mbB_M$ and sample covariance matrix  $\bS_M$ are compared here under various simulation designs.
We set $p = 100$ and $n = 23400$, which represents the case where transactions are recorded per second within one trading day.
We simulate the observations from the following additive model:
\[
\bY_{t_i} = \bX_{t_i} + \bm \varepsilon_{t_i},
\]
in which the log price $(\bX_t)$ follows from the continuous-time process as in \eqref{eq: simulation model}, and the noise values $(\bm\varepsilon_{t_i})$ are drawn independently from $N(0, 0.0002\bI)$. 
The pre-averaging window length $h$ is taken to be  $\lfloor n^{0.55}\rfloor = 252$.
We use $L_i^{(q)}$ to denote the number of transactions for stock $q$ within time stamp $(t_{i-1}, t_i]$, for $q = 1,\cdots, p$ and $i = 1,\cdots, n$.

We designed two transaction schemes as follows.

\begin{figure}[!htb]
\centering
\includegraphics[width= 6cm, height = 6cm]{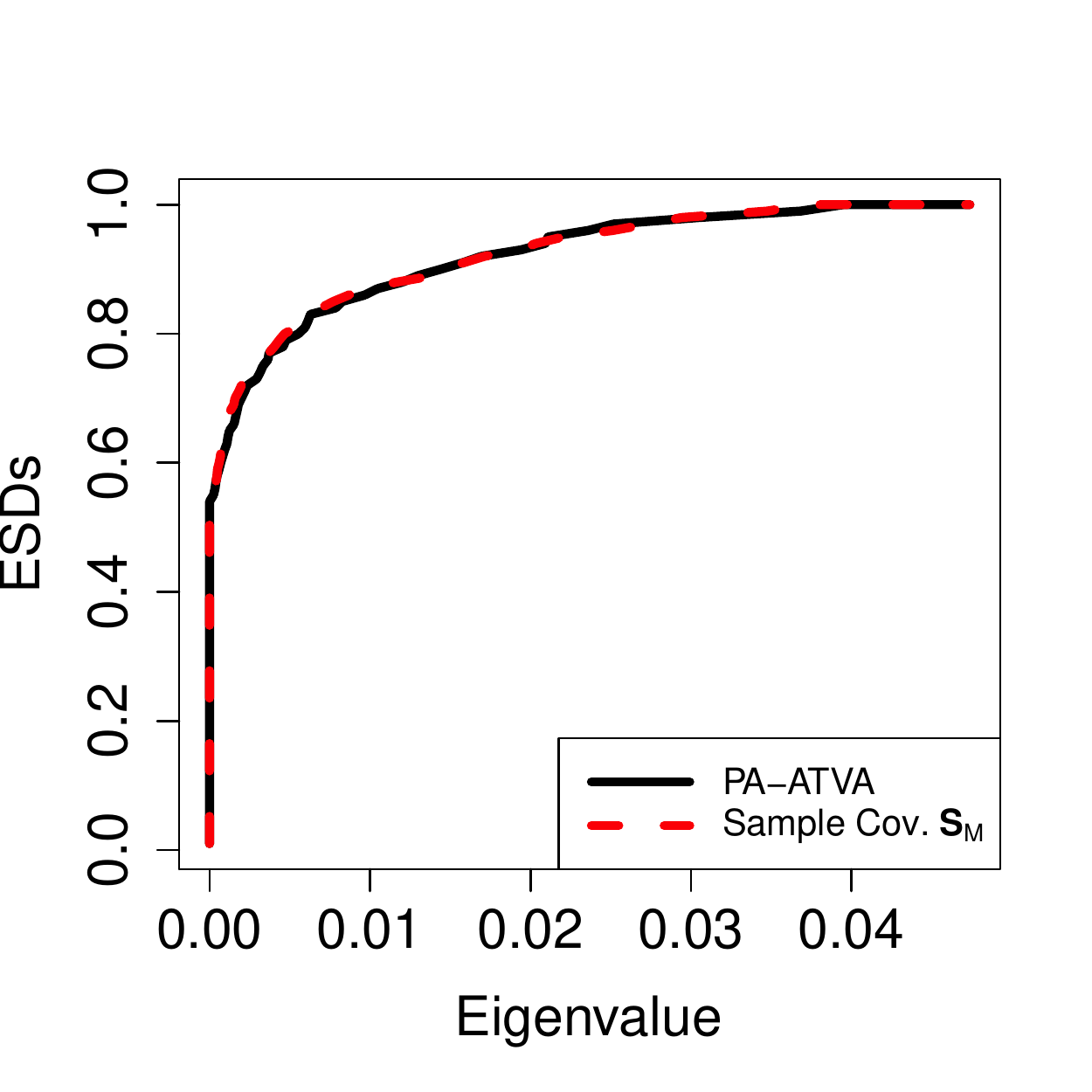}
\caption{
ESDs of PA-ATVA matrix $\mbB_M$ and sample covariance matrix $\bS_M$ for dimension $p  = 100$ and observation frequency $n = 23400$.
The pre-averaging window length $h$ was taken to be $\lfloor n^{0.55}\rfloor = 252$, with an effective sample size $M = \lfloor n/(2h) \rfloor = 46$.
For each stock $q = 1,\cdots, p$, $L_i^{(q)} = L_i$ were generated independently from a Poisson distribution with parameter 20 for the first and the last hours within 6.5 hours of a trading day, and with parameter 5 for the rest of trading hours.
}
\label{fig: application of theorem 2}
\end{figure}

Design \uppercase\expandafter{\romannumeral1}: For simplicity, $L_i^{(q)} = L_i$ for each stock $q$ within time interval $(t_{i-1}, t_i]$, where the $t_i$ values are arranged as an equally spaced grid in $[0,1]$.
The $L_i$'s are generated independently from a Poisson distribution with parameter 20 for the first and the last hours within 6.5 hours of a trading day, and from a Poisson distribution with parameter 5 for the remaining trading hours.
According to the simulation results, shown in Figure \ref{fig: application of theorem 2}, the two ESDs were very closely matched.

\begin{figure}[!htb]
\centering
\includegraphics[width= 6cm, height = 6cm]{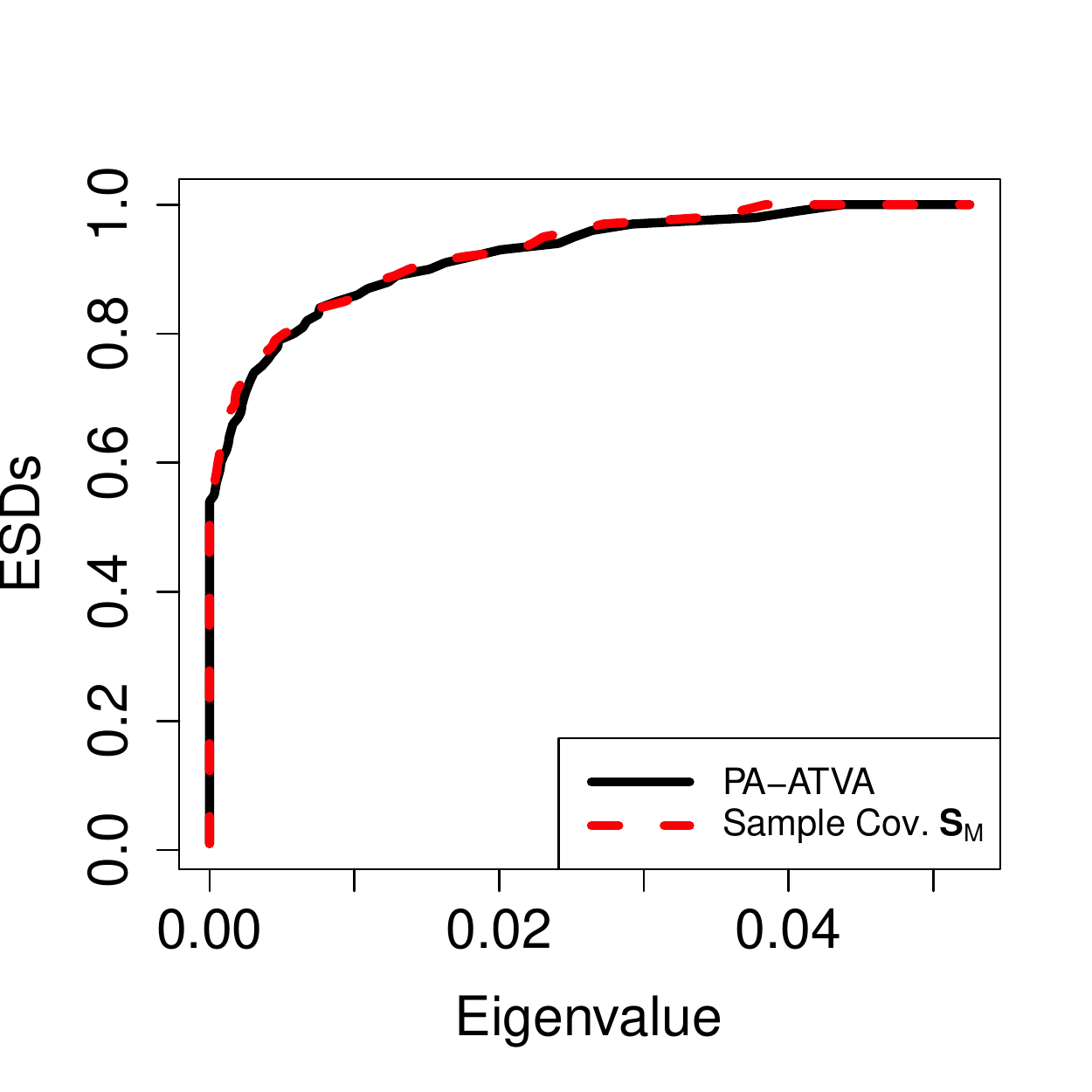}
\caption{ESDs of PA-ATVA matrix $\mbB_M$ and sample covariance matrix $\bS_M$ for dimension $p  = 100$ and observation frequency $n = 23400$.
The pre-averaging window length $h$ was taken to be $\lfloor n^{0.55}\rfloor = 252$, with effective sample size $M = \lfloor n/(2h) \rfloor = 46$.
For each stock $q = 1,\cdots, p$ and each recording time $i = 1,\cdots, n$, $L_i^{(q)}$ values were generated independently from discrete uniform distribution $U[1,L_i]$.
}
\label{fig: application of theorem 2b}
\end{figure}

Design \uppercase\expandafter{\romannumeral2}: In order to generate high-frequency data such as that commonly used in practice, we further simulated observations in a highly asynchronous setting.
Based on Design \uppercase\expandafter{\romannumeral1}, we allowed variation of $L_i^{(q)}$, which is generated independently from a discrete uniform distribution within the interval $[1, L_i]$, for each $q = 1, \cdots, p$ and $i = 1, \cdots, n$.
Figure \ref{fig: application of theorem 2b} displays the ESDs of matrices $\mbB_M$ and $\bS_M$ under Design \uppercase\expandafter{\romannumeral2}.

\subsection{Comparisons with three NS estimators}
In this subsection, we describe some simulation studies to assess the performance of three proposed NS estimators: NS, ANS, and MNS.
We generated a two-day log price process $(\bX_t)$ using model \eqref{eq: simulation model}, except that matrix $\bLambda$ had the following four settings:
\begin{itemize}
  \item Setting \uppercase\expandafter{\romannumeral1}: set $\bLambda = (0.5^{|i - j|})_{i,j = 1,\cdots,p}$ for two days;
  \item Setting \uppercase\expandafter{\romannumeral2}: set $\bLambda = (0.5^{|i - j|})_{i,j = 1,\cdots,p}$ and change the first three eigenvalues of $\bLambda$ to be $15, 10, 5$, and further rescale the spectral of $\bLambda$ to satisfy $\tr(\bLambda\bLambda^{\mT}) = p$;
  \item Setting \uppercase\expandafter{\romannumeral3}: set $\bLambda$ as in Setting \uppercase\expandafter{\romannumeral2}, except that the first three eigenvalues of $\bLambda$ are set to be $15, 10, 5$ for day 1, and the first three eigenvalues of $\bLambda$ are $30, 10, 5$ for day 2.
  \item Setting \uppercase\expandafter{\romannumeral4}:
  set $\bLambda$ as in Setting \uppercase\expandafter{\romannumeral2}, except that the first three eigenvalues of $\bLambda$ are set to be $15, 10, 5$ for day 1, and the first eigenvalue of $\bLambda$ is $30$ for day 2.
\end{itemize}
We constructed high-frequency noisy data $(\bY_{t_i})$ using Design \uppercase\expandafter{\romannumeral2}, described in subsection \ref{subsec: PA-ATVA}, obtaining $23400\times 2$ observations for $p = 30$ or $100$ stocks, respectively.
The pre-averaging window length $h$ in the PA-ATVA matrix $\mbB_M$ and window length $k_n$ in the MNS estimator $\wh{\bSigma}_{\MNS}$ were taken to be  $\lfloor n^{0.55}\rfloor = 252$ and $\lfloor 0.75n^{1/2}\rfloor = 114$, respectively.
For the ANS estimator $\wh{\bSigma}_{\ANS}$, we performed $B = 50$ permutations and compared all the simulation results under the seven data partition criteria given in \eqref{eq: nonlinear shrnkage criterion}.
We denote ANS $i, i = 1,\cdots, 7$, when the data partition criterion $i$ is adopted as in \eqref{eq: critia for selection of n1}.
For any matrix $\mathcal{\bm Q}$, we defined the relative Frobenius loss (RFL) as $||\mathcal{\bm Q}- \ICV||_F/||\ICV||_F$.
For each setting, we repeated the simulations 1000 times. The means and standard deviations for the RFL of these estimators are presented in Table \ref{table: compare the nonlinear shrinkage}.

\begin{table}[!htb]
\small
\centering
\setlength{\abovecaptionskip}{0pt}
\setlength{\belowcaptionskip}{10pt}
\caption{Mean and standard error of the relative Forbenius loss of the different estimators with respect to ICV matrix based on 1000 replications. }
\setlength{\tabcolsep}{1.5 mm}{
\begin{tabular}{cccccccccccc}
  \hline
  \hline
   & \multicolumn{2}{c}{Setting \uppercase\expandafter{\romannumeral1} } && \multicolumn{2}{c}{Setting \uppercase\expandafter{\romannumeral2} } & &\multicolumn{2}{c}{Setting \uppercase\expandafter{\romannumeral3} } &&\multicolumn{2}{c}{Setting \uppercase\expandafter{\romannumeral4}  }\\
   Estimators & mean &std  && mean  &std  && mean  &std && mean  &std\\
     \hline
 \multicolumn{12}{c}{$p = 30, n = 23400$}\\
NS	&0.395	&0.043	&	&0.388	&0.105	&	&0.495	&0.118	&	&0.510	&0.130	\\
ANS2	&0.502	&0.050	&	&0.419	&0.106	&	&0.513	&0.124	&	&0.519	&0.138	\\
ANS3	&0.429	&0.050	&	&0.391	&0.106	&	&0.495	&0.122	&	&0.507	&0.136	\\
ANS4	&0.401	&0.046	&	&0.380	&0.106	&	&0.489	&0.121	&	&0.504	&0.134	\\
ANS5	&0.391	&0.044	&	&0.375	&0.106	&	&0.486	&0.121	&	&0.503	&0.132	\\
ANS6	&0.392	&0.044	&	&0.376	&0.106	&	&0.488	&0.120	&	&0.504	&0.133	\\
ANS7	&0.390	&0.043	&	&0.374	&0.105	&	&0.486	&0.121	&	&0.502	&0.133	\\
MNS	&0.434	&0.034	&	&0.259	&0.065	&	&0.319	&0.126	&	&0.326	&0.137	\\
 \multicolumn{12}{c}{$p = 100, n = 23400$}\\
NS	&0.516	&0.024	&	&0.412	&0.093	&	&0.530	&0.117	&	&0.549	&0.122	\\
ANS1	&0.651	&0.023	&	&0.476	&0.096	&	&0.558	&0.122	&	&0.567	&0.126	\\
ANS2	&0.657	&0.023	&	&0.481	&0.096	&	&0.560	&0.122	&	&0.569	&0.126	\\
ANS3	&0.579	&0.028	&	&0.429	&0.096	&	&0.533	&0.122	&	&0.548	&0.127	\\
ANS4	&0.538	&0.028	&	&0.410	&0.095	&	&0.524	&0.121	&	&0.542	&0.126	\\
ANS5	&0.521	&0.025	&	&0.401	&0.094	&	&0.522	&0.120	&	&0.540	&0.126	\\
ANS6	&0.524	&0.026	&	&0.403	&0.095	&	&0.522	&0.120	&	&0.540	&0.125	\\
ANS7	&0.519	&0.025	&	&0.400	&0.095	&	&0.521	&0.120	&	&0.540	&0.125	\\
MNS	&0.606	&0.022	&	&0.322	&0.053	&	&0.347	&0.105	&	&0.357	&0.119	\\

  \hline
  \hline
\end{tabular}
}
\label{table: compare the nonlinear shrinkage}
\end{table}

As shown in Table \ref{table: compare the nonlinear shrinkage}, NS and ANS7 with small error mean performed better than the other estimators for Setting \uppercase\expandafter{\romannumeral1}, where the underlying covolatility structure was stable and no spikes were considered.
However, when the spiked model of the ICV matrix or structure changes were considered, as in Setting \uppercase\expandafter{\romannumeral2}, \uppercase\expandafter{\romannumeral3}, and \uppercase\expandafter{\romannumeral4}, the MNS estimator outperformed the other estimators in terms of both mean and standard deviation.
Comparing all the ANS estimators with different choices of data splitting criteria, the ANS estimator tended to have better performance when the first data set was larger.
However, generally speaking, the MNS estimator was more efficient and robust when the structure changed or under factor models.

\subsection{Empirical applications}
Given the simulation results described in the previous subsection, we next constructed a minimum variance portfolio using the three proposed NS estimators on intraday high-frequency data for the Dow Jones Industrial Average (DJIA).
We collected tick-by-tick 30 DJIA stock prices from January 1, 2016 to December 31, 2016, comprising 252 trading days, from the Trade and Quote database.

\begin{table}[!htb]
\small
\centering
\begin{threeparttable}
\setlength{\abovecaptionskip}{0pt}
\setlength{\belowcaptionskip}{10pt}
\caption{Summary of the prices for DJIA from January 1, 2016 to December 31, 2016.}
\setlength{\tabcolsep}{1.0 mm}{
\begin{tabular}{lccccc}
  \hline
  \hline
  Corporate	&Total	&Eff. 10-seconds	&Average obs.	&Max. obs.	&per\\
  \hline
3M	&4948363	&526992	&9.39	&548	&85.90\\
American Express	&9355637	&546208	&17.13	&2185	&91.39\\
Apple	&53080066	&588694	&90.17	&8825	&99.96\\
Boeing	&8516963	&553487	&15.39	&1856	&91.51\\
Caterpillar	&9902961	&557979	&17.75	&825	&93.26\\
Chevron	&14862178	&583537	&25.47	&2451	&98.04\\
Cisco Systems	&23185678	&584077	&39.70	&1808	&98.03\\
Coca-Cola	&15540179	&578058	&26.88	&2776	&96.11\\
ExxonMobil	&19374092	&585050	&33.12	&2360	&98.49\\
Goldman Sachs	&9177527	&553721	&16.57	&771	&92.39\\
The Home Depot	&10344101	&575078	&17.99	&650	&95.54\\
IBM	&8224719	&560945	&14.66	&986	&92.32\\
Intel	&25203925	&583081	&43.23	&3945	&97.90\\
Johnson \& Johnson	&13186553	&584239	&22.57	&1506	&97.98\\
JPMorgan Chase	&26160690	&586486	&44.61	&2608	&99.11\\
McDonald's	&9872997	&570740	&17.30	&1725	&94.58\\
Merck	&14945942	&576447	&25.93	&1397	&96.13\\
Microsoft	&39836829	&587821	&67.77	&4651	&99.54\\
Nike	&13927827	&575189	&24.21	&1066	&95.77\\
Pfizer	&27497813	&586176	&46.91	&4897	&98.84\\
Procter \& Gamble	&15151131	&581224	&26.07	&1414	&97.17\\
Travelers	&3934222	&461529	&8.52	&519	&80.50\\
UnitedHealth Group	&7203345	&547476	&13.16	&770	&90.58\\
United Technologies	&8215206	&554497	&14.82	&2247	&91.62\\
Verizon	&18677248	&584807	&31.94	&1298	&98.24\\
Visa	&14610968	&578756	&25.25	&2352	&96.78\\
Wal-Mart	&14691129	&579230	&25.36	&1879	&96.51\\
Walt Disney	&13610611	&583584	&23.32	&1659	&97.80\\
Du Pont	&5753121	&497033	&11.57	&873	&85.28\\
General Electric	&25791517	&586938	&43.94	&2444	&99.14\\

  \hline
  \hline
\end{tabular}
}
\label{table: descriptive statistics of real data}
\begin{tablenotes}
\item[*] ``Total"  means total observations during the selected period. ``Eff. 10-seconds" stands for the number of ten-seconds that contains at least one observations.  ``Averged obs." and ``Max. obs." are the averaged and maximum observed number within the effective ten-seconds. "Per" is the percentage of 10-seconds that contain at least two transactions out of the total effective 10-seconds.
\end{tablenotes}
\end{threeparttable}
\end{table}
Table \ref{table: descriptive statistics of real data} shows a brief summary of the selected dataset.
These stocks had excellent liquidity over the sampling period.
In order to address microstructure noise and asynchronous trading, we adopted the pre-averaging method and set the recording time stamp to 10 seconds.
To clean the data, we filled missing prices using the previous price.
Overnight returns and the first five minutes of trading data were removed for purposes of the calculations, to avoid the effects of price changes and irrelevant jumps.
The minimum variance portfolio solves the following optimization problem
\[
\min\limits_{\bw} \bw^{\mT} (\ICV) \bw, ~~\text{subject to}~\bw^{\mT}\bm{1}_p = 1,
\]
where $\bm{1}_p $ is a vector of $p$ ones and $\ICV$ is a given integrated covolatility matrix.
The solution of $\bw$ is
\begin{equation}\label{eq: the portfolio stratege}
\bw = \frac{(\ICV)^{-1}\bm{1}_p}{\bm{1}_p^{\mT} (\ICV)^{-1}\bm{1}_p}.
\end{equation}
Initially, we invested one unit of capital using \eqref{eq: the portfolio stratege} by constructing three NS estimators, NS, ANS, and MNS.
We re-evaluated portfolio weights every day by considering 10-day (two-week), 15-day (three-week), and 20-day (four-week) training windows and taking the 15-minute pre-averaging window length.
We compared the annualized returns $\hat\mu$, standard deviation $\hat{\sigma}$, and Sharpe ratio $\hat{\mu}/\hat{\sigma}$ for each method.
We also reported the averaged maximum absolute (AME) value of the weight named as maximum exposures. 
For $\ell$-week training windows, $\hat\mu$ and $\hat\sigma$ are defined as
\begin{equation}\label{eq: definition of annalized return}
\hat\mu = \frac{252}{251 - \ell}\sum_{i = \ell + 1}^{251}\bw_i^{\mT}\bm r_i, ~ \hat{\sigma} = \left\{\frac{252}{251 - \ell}\sum_{i = \ell + 1}^{251}(\bw_i^{\mT}\bm r_i - \big(\frac{\hat\mu}{251 - \ell}\big)^2\right\}^{\frac{1}{2}},
\end{equation}
where $\bw_i $ and $\bm r_i$ are the portfolio weights and returns, respectively, for the $i$th day.
For the ANS estimator, we set $B =50$ permutations and further searched the split location $M_{\tau_1}$ by solving the optimization problem \eqref{eq: nonlinear shrnkage criterion}.
For the MNS estimator, we chose $k_n$ to be six minutes to pre-average the transformed data in the estimator $\wh{d_i^{\APA}}$ .
We compared two trading strategies: (1) open-open stands for selling and buying stocks five minutes after the opening of the market; (2) close-close stands for selling and buying stocks at the close of the market.
The empirical results are displayed in Table \ref{table: result for the real data}.
\begin{table}[!htb]
\scriptsize
\centering
\setlength{\abovecaptionskip}{0pt}
\setlength{\belowcaptionskip}{10pt}
\caption{Results of the analysis for Dow Jones Industrial Average stocks.
Return and standard deviation are given in \eqref{eq: definition of annalized return}.
NS, ANS, MNS are three proposed NS estimators in section \ref{subsection: NS} , \ref{subsection: ANS} and \ref{subsection: MNS}, respectively.}
\setlength{\tabcolsep}{1.9 mm}{
\begin{tabular}{llllllllll}
  \hline
  \hline
   Method&Return &Std Dev &Sharpe &AME &&Return &Std Dev &Sharpe &AME \\
   &\%&\%& &\% & &\% &\%& &\%\\
      \hline
   & \multicolumn{4}{c}{open - open} && \multicolumn{4}{c}{close - close}\\
  &\multicolumn{9}{c}{10-day}\\
NS 	&7.36	&10.32	&0.71	&29.01	&	&12.46	&10.51	&1.19	&29.01	\\
ANS	&7.40	&9.46	&0.78	&12.73	&	&9.65	&9.57	&1.01	&12.79	\\
MNS	&10.28	&10.09	&1.02	&18.46	&	&7.29	&10.33	&0.71	&18.46	\\
EW	&19.90	&20.98	&0.95	&3.33	&	&20.03	&14.58	&1.37	&3.33	\\
   &\multicolumn{9}{c}{15-day}\\
NS 	&6.11	&10.20	&0.60	&28.23	&	&8.07	&10.03	&0.80	&28.23	\\
ANS	&11.67	&9.22	&1.27	&13.06	&	&12.98	&9.21	&1.41	&13.07	\\
MNS	&14.46	&10.20	&1.42	&18.37	&	&13.68	&9.67	&1.42	&18.37	\\
EW	&19.20	&20.80	&0.92	&3.33	&	&21.05	&14.48	&1.45	&3.33	\\
      &\multicolumn{9}{c}{20-day}\\
NS 	&7.69	&9.92	&0.77	&26.00	&	&5.25	&9.79	&0.54	&26.00	\\
ANS	&12.86	&9.26	&1.39	&13.12	&	&9.55	&9.00	&1.06	&13.12	\\
MNS	&12.57	&9.75	&1.29	&18.40	&	&9.54	&9.20	&1.04	&18.40	\\
EW	&17.99	&21.00	&0.86	&3.33	&	&17.69	&14.30	&1.24	&3.33	\\
  \hline
  \hline
\end{tabular}
}
\label{table: result for the real data}
\end{table}
The ANS estimator had the smallest standard deviation and smallest maximum exposure for all settings.
With 10-day training windows, the ANS estimator had the largest Sharpe ratio, but on average that of the MNS estimator was larger.
Notably, the ANS and MNS estimators also had much lower computational costs compared with the NS estimator.

\section{Conclusion\label{sec:conclusion} }
This article considers the estimation of high-dimensional ICV matrices based on multiple high-frequency observations.
First, using random matrix theory, we investigated the LSDs between the ATVA matrix and the targeting ICV matrix under the latent log price process and noisy observations, respectively.
Surprisingly, the desirable property of TVA matrix established in \cite{ZhengandLi2011}, that the LSD of TVA matrix depended solely on that of ICV matrix, no longer held.
Our theoretical results show that the LSD of the averaged version of TVA matrix depends not only on that of the ICV matrix but also on the time variability of the number of multiple transactions.
Further, we investigated the limiting spectral property of TVA matrix on the pre-averaging approach and found that this approach worked well, that is, the proposed PA-ATVA matrix eliminated the effects of microstructure noise and asynchronous trading within one recording time stamp.
Therefore, three types of NS estimators were proposed based on the PA-ATVA matrix.
Both simulation and empirical studies indicated that the three estimators (NS, ANS, and MNS) performed reasonably well.
In high-frequency portfolio allocation applications, in particular, the ANS estimator showed the smallest annualized risk and the smallest maximum exposure, but the MNS estimator achieved the largest Sharpe ratio on average.

This work represents a first step in investigating the effects of multiple transactions in a high-dimensional setting.
The results of this paper are encouraging in this direction, and we hope to pursue some extensions in future work.


\newpage
\begin{supplement}[id=suppA]
\stitle{Supplement to ``On the estimation of high-dimensional integrated covariance matrices based on high-frequency data with multiple transactions''}
\sdescription{Due to space constraints, all the proofs are relegated to the supplement.}
\end{supplement}


%


\bigskip

\noindent
Moming Wang: School of Statistics and Management, Shanghai, Shanghai University of Finance and Economics,
777 Guo Ding Road, China, 200433. \\
wangmoming0902@163.com

\bigskip
\noindent
Ningning Xia: School of Statistics and Management, Shanghai, Key Laboratory of Financial Information Technology, Shanghai University of Finance and Economics,
777 Guo Ding Road, China, 200433. \\ xia.ningning@mail.shufe.edu.cn

\bigskip
\noindent
Yong Zhou: Key Laboratory of Advanced Theory and Application in Statistics and Data Science, MOE, and Academy of Statistics and Interdisciplinary Sciences and School of Statistics, East China Normal University, Shanghai, China, 200062. \\ yzhou@amss.ac.cn


\newpage
\setcounter{page}{1}
\begin{supplement}[id=suppA]
  \stitle{Supplement to ``On the estimation of high-dimensional integrated covariance matrices based on high-frequency data with multiple transactions
''}


Owing to space constraints, the proofs of Theorem \ref{thm:LSD_TVA}, \ref{thm: theorem for the TVAPA} and Theorem \ref{thm: theorem for the TVAPA asynt} have been relegated to the Supplementary material.
The proposition/lemma/equation \hbox{etc.} numbers below refer to the numbering in the main article.

\appendix

\renewcommand{\baselinestretch}{1.2}
\setcounter{equation}{0}
\renewcommand{\theequation}{\thesection.\arabic{equation}}



%
%

Throughout this paper, $C$ stands for a constant whose value may change from line to line, and $f_n \sim g_n$ means that $f_n/g_n$ converges to 1.

\section{Proof of Theorem 2.1}

Theorem \ref{thm:LSD_TVA} is a direct consequence of the following two propositions.

\begin{prop}
Under the assumptions of Theorem \ref{thm:LSD_TVA}, the ESD of $\wbSigma$, $F^{\wbSigma}$  converges almost surely to the limit $\widetilde{F}$, which is determined by $\breve{H}$ in that its Stieltjes transform $m_{\widetilde{F}}(z)$ satisfies the following equation:
\begin{equation}m_{\widetilde{F}}(z) = \int_{\tau\in \mathbb{R}} \frac{1}{\tau(1 - \rho(1 + zm_{\widetilde{F}}(z))) - z}d\breve{H}(\tau),~\text{for}~z\in \mathbb{C}^+.
\label{eq: converfence of proposition A1}
\end{equation}
\label{prop: convergence of tilde sigma}
\end{prop}

\no\emph{Proof of Proposition \ref{prop: convergence of tilde sigma}.}  Recall that
\[\Delta\widebar{\bX}_{2i} =  \sum_{j=1}^{L_{2i}}a_{2i,j}\Delta_{2i,j}\bX +   \sum_{j=1}^{L_{2i-1}}b_{2i-1,j}\Delta_{2i-1,j}\bX,\]
where $a_{i,j} = 1-\frac{j-1}{L_i},  b_{i,j} = \frac{j-1}{L_i}$, and $\Delta_{i,j}\bX = \bX_{T_{i-1}+j} - \bX_{T_{i-1}+j - 1}$.
Then we separate $\Delta\widebar{\bX}_{2i}$ into two parts,
\begin{equation}
\Delta\wbX_{2i} = \sqrt{w_i}(\bLambda\bZ_{i} + \bV_i),
\label{eq: expression of delta wbx}
\end{equation}
where
\begin{equation}
w_i = \sum_{j=1}^{L_{2i}}a_{2i,j}^2\int_{s_{T_{2i-1} + j -1}}^{s_{T_{2i-1} + j}}\gamma_t^2dt +  \sum_{j=1}^{L_{2i-1}}b_{2i-1,j}^2\int_{ s_{T_{2i-2} + j -1} }^{ s_{T_{2i-2} + j}}\gamma_t^2dt,
\label{eq: definition of wi}
\end{equation}
\[\bZ_i =  \frac{1}{\sqrt{w_i}}\left[\sum_{j=1}^{L_{2i}}a_{2i,j}\int_{s_{T_{2i-1} + j -1} }^{s_{T_{2i-1} + j} }\gamma_td\bW_t +  \sum_{j=1}^{L_{2i-1}}b_{2i-1,j}\int_{s_{T_{2i-2} + j -1} }^{s_{T_{2i-2} + j} }\gamma_td\bW_t \right],\]
and
\[\bV_i =  \frac{1}{\sqrt{w_i}}\left[\sum_{j=1}^{L_{2i}}a_{2i,j}\int_{s_{T_{2i-1} + j -1} }^{s_{T_{2i-1} + j} }\bmu_tdt +  \sum_{j=1}^{L_{2i-1}}b_{2i-1,j}\int_{s_{T_{2i-2} + j -1} }^{s_{T_{2i-2} + j} }\bmu_tdt  \right].\]
Using the above notation, $\wbSigma$ can be rewritten as
\[\wbSigma = \frac{p}{\hn}\sum_{i=1}^{\hn}\frac{\Delta\widebar{\bX}_{2i}\Delta\wbX_{2i}^{\mathrm{T}}}{|\Delta\wbX_{2i}|^2} = \frac{p}{\hn}\sum_{i=1}^\hn \frac{(\bV_i + \bLambda\bZ_i)(\bV_i + \bLambda\bZ_i)^{\mathrm{T}} }{|\bV_i + \bLambda\bZ_i|^2}.\]
Notice that
\[
\frac{1}{L_{2i}} \sum_{j = 1}^{L_{2i}} a_{2i,j} = \frac{1}{2}(1 + \frac{1}{L_{2i}}), ~~ \frac{1}{L_{2i - 1}}\sum_{j = 1}^{L_{2i - 1}} b_{2i-1, j} = \frac{1}{2}(1 - \frac{1}{L_{2i-1}})
\]
and
\[
\frac{1}{L_{2i}} \sum_{j = 1}^{L_{2i}} a_{2i,j}^2  = \frac{1}{3} + \frac{1}{6L_{2i}^2} + \frac{1}{2L_{2i}}, ~\frac{1}{L_{2i - 1}} \sum_{j = 1}^{L_{2i - 1}} b_{2i - 1,j}^2  = \frac{1}{3} + \frac{1}{6L_{2i - 1}^2} - \frac{1}{2L_{2i-1}},
\]
which are all uniformly bounded as $L_i \geq 1$.
By the boundedness of $|\mu_t^{(l)}|$ and $|\gamma_t|$ from Assumption \eqref{asm:bdd_mu} and \eqref{asm:bdd_gamma}, and the equally spaced transaction durations $\Delta s_{i,j}$ from Assumption \eqref{asm:bdd_duration}, there exists $C>0$ such that $|\bV_i^{(l)}|\leq C/\sqrt{\hn}$ for all $i=1,\cdots,\hn,l=1,\cdots, p$, which indicates that the $|\bV_i|$ are uniformly bounded.
By a similar argument to that used in (3.34) in \cite{ZhengandLi20112}, we have
\begin{equation}
\max\limits_{1\leq i \leq \hn} \big|\frac{1}{p}|\bLambda\bZ_i|^2-1\big| \to 0, \hspace{0.5cm} \text{almost surely}.
\label{eq: max_convergence of tilde z}
\end{equation}
By the boundedness of the $|\bV_i|$, we further have
\begin{equation}
\max\limits_{1\leq i \leq \hn} \big|\frac{1}{p}|\bV_i + \bLambda\bZ_i|^2-1\big| \to 0, \hspace{0.5cm} \text{almost surely}.
\label{eq: order of v pluz z}
\end{equation}
By the same arguments used in the proof of Theorem 2 in \cite{ZhengandLi20112}, we can show that Proposition \ref{prop: convergence of tilde sigma} holds.\proofend

\begin{prop}
Under the assumptions of Theorem \ref{thm:LSD_TVA},
the ESD of the ICV matrix converges almost surely in distribution to a probability distribution $H$ as $p\to\infty$ defined by
\begin{equation}
H(x) = \breve{H}(x/\theta),
\label{eq: the relationship of ICV and H}
\end{equation}
where $\theta = \int_0^1(\gamma_t^*)^2dt$.
Moreover,
\begin{equation}
\lim\limits_{p\to\infty} \frac{1}{p}\sum_{i = 1}^N |\Delta\wbX_{2i}|^2 =  \int_0^1 (\dfrac 13 + \dfrac 16 f_{2,s}^*)(\gamma_s^*)^2ds,\hspace{0.5cm} \text{almost surely}.
\label{eq: convergence of Vi plus lambda Zi}
\end{equation}
\label{prop: convergence of wi}
\end{prop}
\no\emph{Proof of Proposition \ref{prop: convergence of wi}.} The convergence of $F^{ICV}$ follows from Assumption \eqref{asm:LSD_bH} and the fact that
\[F^{\ICV}(x) = F^{\breve{\bSigma}}\left(\frac{x}{\int_0^1\gamma_t^2dt}\right) \hspace{0.5cm} \text{for all~} x\geq 0.\]
We proceed to show the convergence of $1/p\cdot \sum_{i=1}^N|\Delta \wbX_{2i}|^2$.
Write $\Delta \wbX_{2i} = \sqrt{w_i}(\bLambda \bZ_i + \bV_i)$ as in \eqref{eq: expression of delta wbx}; then
\[
\dfrac 1p |\Delta\wbX_{2i}|^2 = \dfrac 1p \sum_{i=1}^N w_i |\bV_i + \bLambda \bZ_i|^2 =\sum_{i=1}^N Ew_i + \sum_{i=1}^N (w_i - Ew_i) + \varepsilon,
\]
where $\varepsilon = \sum_{i = 1}^N w_i (p^{-1} |\bV_i + \bLambda \bZ_i|^2 - 1)$.
First, almost surely, that $\max\limits_{i,n} (N w_i)$ is bounded follows from Assumption \eqref{asm:bdd_gamma} and \eqref{asm:bdd_duration}.
Thus, combining this with the convergence of \eqref{eq: order of v pluz z}, the error term $\varepsilon$ converges to 0 almost surely.
Besides, $\sum_{i=1}^N (w_i - E(w_i))$ converges to zero almost surely by the Kolmogorov strong law of large numbers and the boundedness of $\max\limits_{i,n} (N w_i)$.
Next, we show that there exists a piecewise continuous process $(w_s)$ with finite jumps such that
\[
\lim_{n\to\infty}\sum_{i=1}^N \int_{(2i-2)/n}^{2i/n}|N E(w_i) - w_s|ds =0,
\]
where $w_s = (\dfrac 13 + \dfrac 16 f_{2,s}^*)\cdot (\gamma_s^*)^2$.
Suppose that $(\gamma_t^*)$, $(f_{1,s}^*)$, and $(f_{2,s}^*)$ have a finite number of jumps J in total. 
For each $j = 1,\cdots, J$, there exists an $\ell_j$ such that the $j$th jump falls in the interval $((2\ell_j - 2)/n, 2\ell_j/n]$.
Then
\begin{equation*}
\begin{split}
&\sum_{i=1}^\hn\int_{(2i-2)/n}^{2i/n}|N E(w_i) - w_s|ds\\
=&\sum\limits_{\ell_j \in \{\ell_1,\cdots,\ell_J\}}\int_{\frac{2\ell_j-2}{n}}^{\frac{2\ell_j}{n}}|N E(w_{\ell_j}) - w_s|ds +
\sum\limits_{i \notin \{\ell_1,\cdots,\ell_J\}}\int_{\frac{2i-2}{n}}^{\frac{2i}{n}}|N E(w_i) - w_s|ds\\
:= &D_1 + D_2.\\
\end{split}
\end{equation*}
As $(N w_{\ell_j})$, $|\gamma_s^*|$, and $|f_{2,s}^*|$ are bounded, for any $\epsilon >0$ and for sufficiently large $n$, we have
\[|D_1|\leq \frac{2}{n}JC <\epsilon.\]
Notice that
\[
E \left(\frac{1}{L_{2i}}\sum_{j = 1}^{L_{2i}}a_{2i,j}^2 \right)= \frac{1}{3} + \frac{1}{6}E\big(\frac{1}{L_{2i}^2}\big) + \frac{1}{2}E\big(\frac{1}{L_{2i}}\big),
\]
and
\[ E\left( \frac{1}{L_{2i-1}}\sum_{j = 1}^{L_{2i-1}}b_{2i-1
,j}^2 \right)  = \frac{1}{3} + \frac{1}{6} E\big( \frac{1}{L_{2i-1}^2}\big) - \frac{1}{2} E\big( \frac{1}{L_{2i-1}}\big).
\]
For the second term $D_2$, because $(\gamma_t^*)$, $f_{1,s}^*$,  and $f_{2,s}^*$ are continuous in $[(2i-2)/n, 2i/n]$ when $i \notin \{\ell_1,\cdots,\ell_J\}$, and by Assumption \eqref{asm:bdd_ICV} and \eqref{asm:moment_L}, $(\gamma_t)$, $(f_{1,s})$,  and $(f_{2,s})$ uniformly converge to $(\gamma_t^*)$, $(f_{1,s}^*)$,  and $(f_{2,s}^*)$, respectively.
For any $\epsilon >0$ and for sufficiently large $n,p$, we have
\[|\gamma_t^* - \gamma_{2i/n}|<\epsilon, ~~ \text{and}~~ |f_{2, 2i/n}^* - f_{2, t}^*|<\epsilon, ~~ \text{for all}  ~~t \in [\frac{2i-2}{n},\frac{2i}{n} ],\]
for $i = 1,\cdots, N$ and $|\gamma_t - \gamma_t^*|<\epsilon$ for all $t\in [0,1]$.
Define
\[
\begin{split}
w_i^* &= \sum_{j=1}^{L_{2i}}a_{2i,j}^2\int_{s_{T_{2i-1} + j -1}}^{s_{T_{2i-1} + j}}(\gamma_t^*)^2dt +  \sum_{j=1}^{L_{2i-1}}b_{2i-1,j}^2\int_{ s_{T_{2i-2} + j -1} }^{ s_{T_{2i-2} + j}}(\gamma_t^*)^2dt,\\
\tw_i^*&= (\gamma_{2i/n}^*)^2 \frac{1}{nL_{2i}}\sum_{j=1}^{L_{2i}}a_{2i,j}^2 + (\gamma_{2i/n}^*)^2 \frac{1}{nL_{2i-1}}\sum_{j=1}^{L_{2i-1}}b_{2i-1,j}^2\\
\text{and}~~~
 \ttw_i^*&=(\gamma_{2i/n}^*)^2\cdot\frac{2}{n}\bigg(\frac{1}{3} + \frac{1}{6}f_{2, 2i/n}^*\bigg).
\end{split}
\]
For all large $n$, because $|\gamma_t|\leq C$ and $1\leq L_i< L^*$, for any $i, t$, we have
\[
\begin{split}
&|D_2| \\
\leq& \sum\limits_{i}\int_{\frac{2i-2}{n}}^{\frac{2i}{n}}N|E(w_i) - E(w_i^*)|ds + \sum\limits_{i\notin \{\ell_1,\cdots,\ell_J\}}\int_{\frac{2i-2}{n}}^{\frac{2i}{n}}N| E(w_i^*) - E(\tw_i^*)|ds \\
+&\sum\limits_{i}\int_{\frac{2i-2}{n}}^{\frac{2i}{n}}| NE(\tw_i^*) - N\ttw_i^*|ds  +\sum\limits_{i\notin \{\ell_1,\cdots,\ell_J\}}\int_{\frac{2i-2}{n}}^{\frac{2i}{n}}| N\ttw_i^* - w_s|ds\\
\to& 0, ~~~\text{almost surely.}
\end{split}
\]
This completes the proof of \eqref{eq: convergence of Vi plus lambda Zi}.\proofend

\section{Proof of Theorem 2.2}
\label{app: appendix of theorem 2}
The convergence of the ESD of the ICV was proved in Theorem \ref{thm:LSD_TVA}.
The rest of Theorem \ref{thm: theorem for the TVAPA} is a direct consequence of the following two lemmas.
\begin{lem}\label{lem:LSD ttSigma}
Under the assumptions of Theorem \ref{thm: theorem for the TVAPA}, the ESD of $\hbSigma$ converges almost surely, and the limit $F^{\hbSigma}$ is determined by  $\breve{H}$ in that its Stieltjes transform $m_{\hbSigma}(z)$ satisfies the following equation:
\begin{equation}m_{\hbSigma}(z) = \int_{\tau\in \mathbb{R}} \frac{1}{\tau(1 - c(1 + zm_{\hbSigma}(z))) - z}d\breve{H}(\tau),\hspace{0.2cm}\text{for}~z\in \mathbb{C^+}.
\label{eq:LSD of ttF}
\end{equation}

\end{lem}

\begin{lem}\label{lem:theta estimation}
Under the assumptions of Theorem \ref{thm: theorem for the TVAPA}, we have
\begin{equation}\lim_{p\to\infty}3\frac{\sum_{i=1}^\hv|\Delta\wwbY_{2i}|^2}{p} =  \theta, \hspace{0.5cm}\text{almost surely}.
\label{eq: convergence of k avergence of wwX norm}
\end{equation}
\end{lem}

\no\emph{Proof of Lemma \ref{lem:LSD ttSigma}.}
To prove the convergence of $F^{\hbSigma}$, we shall show that
\[
\hbSigma = \frac{p}{\hv} \sum_{i=1}^\hv\frac{\Delta\twbY_{2i}(\Delta\twbY_{2i})^{\mathrm{T}}}{|\Delta\twbY_{2i}|^2}~~~\text{and}~~~
\ttbSigma := \frac{p}{\hv}\sum_{i=1}^{\hv}\frac{\Delta\twbX_{2i}\Delta\twbX_{2i}^{\mathrm{T}}}{|\Delta\twbX_{2i}|^2}
\]
have the same LSD.
Following the same arguments as in the proof of Proposition C.1 of \cite{XiaandZheng20182}, it suffices to show that
\begin{equation}
\max\limits_{1\leq i \leq \hv, 1\leq j \leq p} \frac{\sqrt{p}|\Delta\twbe_{2i}^{(j)}|}{|\Delta\twbX_{2i}|}\to 0,\hspace{1cm} \text{almost surely}.
\label{eq: convergence og delta epsilon divided by delta x}
\end{equation}
We start by showing that there exists a constant $\widetilde{C}>0$ for large $M$, such that
\begin{equation}\label{eq: minimun of twbX}
\min\limits_{1\leq i \leq M} |\Delta\twbX_{2i}|^2 \geq \widetilde{C}.
\end{equation}

By equation \eqref{eqn:Delta_V}, we decompose $\Delta\twbX_{2i}$ as follows:
\begin{align*}
\Delta\twbX_{2i}= & \frac{1}{h}\sum_{j=1}^{h}\left(\wbX_{(2i-1)h + j} - \wbX_{(2i-2)h}\right) - \frac{1}{h}\sum_{j=1}^{h}\left( \wbX_{(2i-2)h+j} -  \wbX_{(2i-2)h}\right) \notag\\
= & \sum_{l=1}^{2h}\left (1 - \frac{|h-l+1|}{h} \right)\Delta\wbX_{(2i-2)h+l}\notag\\
= & \sum_{l=1}^{2h}\left(1 - \frac{|h-l+1|}{h} \right)\sum_{j=1}^{L_{(2i-2)h+l}}a_{(2i-2)h+l,j}\Delta_{(2i-2)h+l,j}\bX\notag\\
&+\sum_{l=1}^{2h-1}\left(1 - \frac{|h-l|}{h} \right )\sum_{j=1}^{L_{(2i-2)h+l}}b_{(2i-2)h+l,j}\Delta_{(2i-2)h+l,j}\bX.\\
\end{align*}
From the fact that  $a_{i,j} + b_{i,j} = 1$, we can further write $ \Delta\twbX_{2i}$ as

\begin{align}\label{eq: expression for tbwX}
\Delta\twbX_{2i} = & \frac{1}{h}\sum_{j=1}^{L_{2ih}}a_{2ih,j}\Delta_{2ih,j}\bX \notag\\
&+ \sum_{l=1}^{2h-1}\sum_{j=1}^{L_{(2i-2)h + l}} \left(1-\frac{|h-l+1|}{h}\right)\Delta_{(2i-2)h+l,j}\bX \notag\\
&+ \sum_{l=1}^{2h-1}\sum_{j=1}^{L_{(2i-2)h + l}} \left[ \frac{|h-l+1| - |h-l|}{h} b_{(2i-2)h+l,j}\right]\Delta_{(2i-2)h+l,j}\bX \notag\\
:= & \tbV_i + \sqrt{\psi_i}\bLambda\tbZ_i,
\end{align}
where
\begin{align*}
\psi_i &:= \frac{1}{h^2}\sum_{j=1}^{L_{2ih}}a_{2ih,j}^2\int_{s_{T_{2ih-1} + j - 1}}^{s_{T_{2ih-1} + j}}\gamma_t^2dt\\
&+\sum_{l=1}^{2h-1}\sum_{j=1}^{L_{(2i-2)h + l}} \left[\left(1-\frac{|h-l+1|}{h}\right) + \frac{|h-l+1| - |h-l|}{h} b_{(2i-2)h+l,j}\right]^2\\
&\hspace{2cm}\cdot\int_{s_{T_{(2i-2)h+l-1} + j - 1}}^{s_{T_{(2i-2)h+l-1} + j}} \gamma_t^2dt,
\end{align*}
\begin{align*}
\tbZ_i& := \frac{1}{\sqrt{\psi_i}}\frac{1}{h}\sum_{j=1}^{L_{2ih}}a_{2ih,j}\int_{s_{T_{2ih-1} + j - 1}}^{s_{T_{2ih-1} + j}}\gamma_td\bW_t\notag\\
&+\sum_{l=1}^{2h-1}\sum_{j=1}^{L_{(2i-2)h + l}} \left[\left(1-\frac{|h-l+1|}{h}\right) + \frac{|h-l+1| - |h-l|}{h} b_{(2i-2)h+l,j}\right] \\
&\hspace{2cm}\cdot\int_{s_{T_{(2i-2)h+l-1} + j - 1}}^{s_{T_{(2i-2)h+l-1} + j}} \gamma_td\bW_t,
\end{align*}
and
\begin{align*}
\tbV_i &:= \frac{1}{h}\sum_{j=1}^{L_{2ih}}a_{2ih,j}\int_{s_{T_{2ih-1} + j - 1}}^{s_{T_{2ih-1} + j}}\bmu_tdt\notag\\
&+\sum_{l=1}^{2h-1}\sum_{j=1}^{L_{(2i-2)h + l}} \left[\left(1-\frac{|h-l+1|}{h}\right) + \frac{|h-l+1| - |h-l|}{h} b_{(2i-2)h+l,j}\right] \\
&\hspace{2cm}\cdot\int_{s_{T_{(2i-2)h+l-1} + j - 1}}^{s_{T_{(2i-2)h+l-1} + j}}\bmu_tdt.
\end{align*}
Without loss of generality, we may assume that $\gamma_t$ and $\bW_t$ are independent, which leads to the fact that each entry of $\tbZ_i$ is i.i.d. standard normal.
Otherwise, by using a similar trick as in the proof of (3.34) of \cite{ZhengandLi20112}, we have
\begin{equation}\label{eq: norm convergeb=nce of tbZ}
\max\limits_{1\leq i \leq M} |\frac{1}{p}|\bLambda\tbZ_i|^2  - 1 | \to 0, ~~~\text{almost surely}.
\end{equation}
Combining this with the fact that all the entries of $\tbV_i$ are of order $O(h/n) = o(1/\sqrt{p})$, we have
\begin{equation}
\frac{\sum_{i=1}^\hv|\Delta\twbX_{2i}|^2}{p} = \frac{\sum_{i=1}^\hv|\tbV_i + \sqrt{\psi_i}\bLambda\tbZ_i|^2}{p} = \sum_{i=1}^\hv\psi_i + o_{a.s.}(1).
\label{eq: psii plus error term}
\end{equation}

Thus, from equation \eqref{eq: expression for tbwX},  we have
\[
|\Delta\twbX_{2i}|^2 = |\tbV_i + \sqrt{\psi_i}\bLambda\tbZ_i|^2 \geq |\tbV_i|^2 + |\psi_i|\cdot| \bLambda\tbZ_i|^2 - 2|\tbV_i| | \sqrt{\psi_i}\bLambda\tbZ_i|.
\]
Assumption \eqref{asm:bdd_sigma} implies that for all $i$, there exists $C'$ such that$|\psi_i|\geq C'h/n.$
Taking this together with Assumption \eqref{asm:dim_windows} and Equation \eqref{eq: norm convergeb=nce of tbZ}, there exists $C^*>0$ such that for all $n$ large enough,
\[
\min\limits_{1\leq i \leq m}|\psi_i|\cdot| \bLambda\tbZ_i|^2 \geq C^*.
\]
Moreover, $\max_i |\tbV_i| =O(\sqrt{p}\times h/n) = o(1)$ follows from Assumption \eqref{asm:bdd_mu}.
Therefore \eqref{eq: minimun of twbX} follows.
Next, we will show that
\begin{equation}\label{eq: sqrtt p delta twbe}
\max\limits_{1\leq i \leq M, 1\leq j \leq p} \sqrt{p}|\Delta \twbe_{2i}^{(j)} | \to 0, ~~~ \text{almost surely}.
\end{equation}
Observe that $(\bar{\varepsilon}_i^{(j)}) = (1/L_i\sum_{k=1}^{L_i}\varepsilon_{i,k}^{(j)})$ is also a $\rho$-mixing sequence following from Definition \ref{def: rho-misiing}, and the $\rho$-mixing coefficients $\bar{\rho}^j(r)$ based on $(\bar{\varepsilon}_i^{(j)})$ have the same order as $\rho^j(r)$ by the boundedness of $L_i$ from Assumption \eqref{asm:bdd_L}.
Therefore, \eqref{eq: sqrtt p delta twbe} follows by the same process used in (C.7) of \cite{XiaandZheng20182}, which, together with \eqref{eq: minimun of twbX}, implies that \eqref{eq: convergence og delta epsilon divided by delta x} holds.\proofend

\no\emph{Proof of Lemma \ref{lem:theta estimation}.}
Note that
\[\sum_{i=1}^\hv|\Delta\twbY_{2i}|^2 = \sum_{i=1}^\hv|\Delta\twbX_{2i}|^2 + 2\sum_{i=1}^\hv\Delta\twbX_{2i}^{\mathrm{T}}\Delta\twbe_{2i} + \sum_{i=1}^M|\Delta\twbe_{2i}|^2.
\]
The convergence of \eqref{eq: convergence og delta epsilon divided by delta x} implies that $\sum_{i=1}^M|\Delta \twbe_{2i}|^2/p \to 0$ almost surely.
It remains to prove that
\begin{equation}\label{eq: convergence of sum of TWBx}
\lim_{p\to\infty} 3\frac{\sum_{i=1}^\hv|\Delta\twbX_{2i}|^2}{p} = \theta, \hspace{1cm} \text{almost surely},
\end{equation}
and
\begin{equation}\label{eq: twbX plus twbe}
\frac{\sum_{i=1}^\hv\Delta\twbX_{2i}^{\mathrm{T}}\Delta\twbe_{2i}}{p}\to 0, \hspace{1cm}\text{almost surely}.
\end{equation}
To show \eqref{eq: convergence of sum of TWBx}, by \eqref{eq: psii plus error term}, it suffices to show that
\[
\lim_{n\to\infty}\sum_{i=1}^\hv\int_{s_{T_{(2i-2)h}}}^{s_{T_{2ih}}}|\hv\psi_i - \frac{1}{3}(\gamma_s^*)^2|ds =0,
\]
almost surely.
Suppose that  $\gamma_t^*$ has $J$ jumps at $\{\tau_1,\cdots,\tau_J\}$; then
\begin{align*}
&\sum_{i=1}^\hv\int_{s_{T_{(2i-2)h}}}^{s_{T_{2ih}}} |\hv\psi_i - \frac{1}{3}(\gamma_s^*)^2|ds\\
=&\sum\limits_{i \in \{\tau_1,\cdots,\tau_J\}}\int_{s_{T_{(2i-2)h}}}^{s_{T_{2ih}}}|\hv\psi_i - \frac{1}{3}(\gamma_s^*)^2|ds \\
&+ \sum\limits_{i \notin \{\tau_1,\cdots,\tau_J\}}\int_{s_{T_{(2i-2)h}}}^{s_{T_{2ih}}}|\hv\psi_i - \frac{1}{3}(\gamma_s^*)^2|ds \\
:= &\Delta_1 + \Delta_2.
\end{align*}
For any $\epsilon > 0$ and for sufficiently large $n$, $|\Delta_1| \leq \epsilon$ follows from the boundedness of $|M\psi_i|$ and $\gamma_t^*$.
For the second term, $\Delta_2$, by defining $\psi_i^*$ by replacing $\gamma_t$ with $\gamma_t^*$ from the definition of $\psi_i$, we have
\[
|\Delta_2| \leq \Delta_{21} + \Delta_{22} + \Delta_{23} +\Delta_{24},
\]
where
\[
\begin{split}
\Delta_{21} &:= \sum\limits_{i \notin \{\tau_1, \cdots, \tau_J\}} \int_{s_{T_{(2i-2)h}}}^{s_{T_{2ih}}} |M\psi_i - M\psi_i^*|ds,\\
\Delta_{22} &:= \sum\limits_{i \notin \{\tau_1, \cdots, \tau_J\}} \int_{s_{T_{(2i-2)h}}}^{s_{T_{2ih}}} | M\psi_i^* - M(\gamma_{(2i-2)/h}^*)^2A_i|ds,\\
\Delta_{23} &:= \sum\limits_{i \notin \{\tau_1, \cdots, \tau_J\}} \int_{s_{T_{(2i-2)h}}}^{s_{T_{2ih}}} | M(\gamma_{(2i-2)/h}^*)^2A_i - M(\gamma_{s}^*)^2A_i |ds,\\
\Delta_{24} &:= \sum\limits_{i \notin \{\tau_1, \cdots, \tau_J\}} \int_{s_{T_{(2i-2)h}}}^{s_{T_{2ih}}} | M(\gamma_{s}^*)^2A_i - \frac{1}{3}(\gamma_{s}^*)^2 |ds,
\end{split}
\]
and
\[
\begin{split}
A_i := &\frac{1}{h^2} \sum_{j=1}^{L_{2ih}} a_{2ih,j}^2 \Delta s_{2ih,j} + \sum_{l=1}^{2h-1} \sum_{j=1}^{L_{(2i-2)h+l}} \Big[1 - \frac{|h-l+1|}{h}+ \\
&\hspace{1cm}+\frac{|h-l+1| - |h-l|}{h}b_{(2i-2)h+l,j} \Big]^2\Delta s_{(2i-2)h+l,j} .
\end{split}
\]
We further decompose $A_i$ as $A_i = A_{i1} +  A_{i2} + A_{i3}+ A_{i4}$, where
\begin{align*}
A_{i1} &:= \frac{1}{h^2} \sum_{j=1}^{L_{2ih}} a_{2ih,j}^2 \Delta s_{2ih,j} = O(\frac{1}{nh^2}),\\
A_{i2} &:= \sum_{l=1}^{2h-1} \left(1 - \frac{|h-l+1|}{h} \right)^2\sum_{j=1}^{L_{(2i-2)h+l}} \Delta s_{(2i-2)h+l,j}= \frac{2h}{3n} + o(\frac{h}{n}),\\
A_{i3} &:= \sum_{l=1}^{2h-1} \sum_{j=1}^{L_{(2i-2)h+l}} \frac{2}{h}\left(1 - \frac{|h-l+1|}{h} \right)(|h-l+1| - |h-l|) \\
&\cdot b_{(2i-2)h+l,j}\Delta s_{(2i-2)h+l,j} = O(\frac{1}{n}),\\
A_{i4} &:= \sum_{l=1}^{2h-1} \sum_{j=1}^{L_{(2i-2)h+l}} \frac{1}{h^2}(|h-l+1| - |h-l|)^2b_{(2i-2)h+l,j}^2\Delta s_{(2i-2)h+l,j} \\
&= O(\frac{1}{nh}),
\end{align*}
follow from the boundedness of $L_i$ both below and above, and $\Delta s_{i,j} = O(1/n)$ by Assumption \eqref{asm:interval}.
As $(\gamma_t^*)$ is continuous in $[s_{T_{(2i-2)h}}, s_{T_{2ih}}]$ when $i \notin \{\tau_1, \cdots, \tau_J\}$, $(\gamma_t)$ uniformly converges to $(\gamma_t^*)$  by Assumption \eqref{asm:bdd_sigma}, and $\psi_i^* = O(h/n)$ by Assumption \eqref{asm:bdd_gamma} and \eqref{asm:bdd_L}, for any $\epsilon >0 $ and sufficiently large $n,p$, it is easy to show that
\[
|MA_i - 1/3|\leq \epsilon ~~~\text{and}~~~\max\{\Delta_{21},\Delta_{22}, \Delta_{23},\Delta_{24}  \}<C\epsilon.
\]
This completes the proof of \eqref{eq: convergence of sum of TWBx}.
Finally, \eqref{eq: twbX plus twbe} follows from \eqref{eq: convergence of sum of TWBx} and \eqref{eq: convergence og delta epsilon divided by delta x}.\proofend

\section{Proof of Theorem 2.3}
Without of generality, we may assume that  $\gamma_t$ is independent of $\bW_t$.
Suppose that we have $L_i^{(q)}(\geq 1)$ observations for each stock $q$ at recording time $t_i = i/n$, for $q = 1,2,\cdots, p$ and $i = 1,2,\cdots, n$.
Recall that for any process $(\bV_t)$, $V_{i,j}^{(q)}$ denote the observation of the $j$th transaction for stock $q$ during time interval $(t_{i-1}, t_i]$,
and the true transaction time of $V_{i,j}^{(q)}$ is denoted as $s_{T_{i-1}+j}^{(q)}$, for $j = 1,\cdots, L_i^{(q)}$ satisfying $t_{i-1} \leq s_{T_{i-1}}^{(q)}< s_{T_{i-1} + 1}^{(q)}<\cdots < s_{T_{i-1} + L_i^{(q)}}^{(q)} = s_{T_i}^{(q)} \leq t_i$.
Under asynchronous trading conditions, the average of multiple observations at each recording time $t_i$ is denoted by
\[
\wbV^*_i := \left(\sum_{j=1}^{L_i^{(1)}}\frac{1}{L_i^{(1)}}V_{i,j}^{(1)}, \cdots, \sum_{j=1}^{L_i^{(p)}}\frac{1}{L_i^{(p)}}V_{i,j}^{(p)} \right)^{\mT};
\]
thus, the increment at trading time $t_i$ becomes
\[\Delta\widebar{\bV}^*_i :=  \begin{pmatrix}
\sum_{j=1}^{L_i^{(1)}}a_{i,j}^{(1)}\Delta_{i,j}^{(1)}V +   \sum_{j=1}^{L_{i-1}^{(1)}}b_{i-1,j}^{(1)}\Delta_{i-1,j}^{(1)}V\\
\vdots\\
\sum_{j=1}^{L_i^{(p)}}a_{i,j}^{(p)}\Delta_{i,j}^{(p)}V +   \sum_{j=1}^{L_{i-1}^{(p)}}b_{i-1,j}^{(p)}\Delta_{i-1,j}^{(p)}V\\
\end{pmatrix},\]
where $a_{i,j}^{(q)} = 1 - \frac{j-1}{L_i^{(q)}}$, $b_{i,j}^{(q)} =  \frac{j-1}{L_i^{(q)}}$ and  $\Delta_{i,j}^{(q)}X = X_{i, j}^{(q)} -  X_{i, j-1}^{(q)}$, for $q = 1,2,\cdots,p$.
Similar to the proof of Lemma \ref{lem:LSD ttSigma}, we first show that
there exists a constant $C>0$ such that for large $M$
\begin{equation}\label{eq: xy same require in AT}
\min\limits_{1\leq i\leq M} |\Delta \twbX^*_{2i}|^2 \geq C,
\end{equation}
and
\begin{equation}\label{eq: xy same require in AT done}
\max\limits_{1\leq i \leq M, 1\leq j \leq p} \sqrt{p}|\Delta (\twbe_{2i}^{*})^{(j)} | \to 0, ~~~\text{almost surely,}
\end{equation}
which leads to the result that the PA-ATVA matrix based on noisy observations $\Delta\twbY^*_{2i}$ and that based on latent log price $\Delta\twbX^*_{2i}$ have the same LSD.
We only need to prove \eqref{eq: xy same require in AT}, as \eqref{eq: xy same require in AT done} holds straightforwardly from the proof of \eqref{eq: sqrtt p delta twbe} given earlier.
By a similar decomposition to that used in \eqref{eq: expression for tbwX}, the return based on the pre-averaged (latent) price can be decomposed as
\[
\Delta\twbX_{2i}^* = \frac{1}{h}\sum_{j=1}^{h}\wbX_{(2i-1)h+j}^* - \frac{1}{h}\sum_{j=1}^{h}\wbX_{(2i-2)h+j}^* = \bR_{i1} + \bM_{i} + \bR_{i2},
\]
where $\bM_i = (M_{i}^{(1)}, \cdots, M_{i}^{(q)})^{\mT}$, $\bR_{i\ell} = (R_{i\ell}^{(1)}, \cdots, R_{i\ell}^{(q)})^{\mT}$ for $\ell = 1, 2$, and their $q$th components have the form
\begin{align*}
R_{i1}^{(q)} = &\frac{1}{h}\sum_{j=1}^{L_{2ih}^{(q)}}a_{2ih,j}^{(q)}\Delta_{2ih,j}^{(q)}X,\\
M_{i}^{(q)} = & \sum_{l=1}^{2h-1} \sum_{j = 1}^{L_{(2i-2)h+l}^{(q)} }c_{l,h} \Delta_{(2i-2)h + l,j}^{(q)}X,\\
R_{i2}^{(q)} = &\sum_{l=1}^{2h-1}\sum_{j=1}^{L_{(2i-2)h + l}^{(q)}} \left[ \frac{|h-l+1| - |h-l|}{h}
b_{(2i-2)h+l,j}^{(q)}\right]\Delta_{(2i-2)h+l,j}^{(q)} X,
\end{align*}
where $c_{l,h} = 1 - \frac{|h-l+1|}{h} $.
Note that $M_i^{(q)}$ can be reduced to
\[
\sum_{l=1}^{2h-1} c_{l,h} \cdot (X_{s_{T_{(2i-2)h + l} }^{(q)}}^{(q)} - X_{s_{T_{(2i-2)h + l - 1} }^{(q)}}^{(q)}).
\]
We further decompose $M_i^{(q)}$ as $W_i^{(q)} + R_{i3}^{(q)}$, where
\begin{align*}
W_{i}^{(q)} = &\sum_{l=1}^{2h-1} c_{l,h} (X_{t_{(2i-2)h + l}}^{(q)} -  X_{t_{(2i-2)h + l - 1}}^{(q)}),\\
R_{i3}^{(q)} = &\sum_{l=1}^{2h-1} c_{l,h} \left(  X_{s_{T_{(2i-2)h + l} }^{(q)}}^{(q)} - X_{t_{(2i-2)h + l} }^{(q)} + X_{t_{(2i-2)h + l -1}}^{(q)} -  X_{s_{T_{(2i-2)h + l - 1}}^{(q)}}^{(q)}  \right),
\end{align*}
and $X_{t_i}^{(q)}$ denotes the log price for stock $q$ at recording time $t_i$.  Let $\bW_i = (W_{i}^{(1)}, \cdots, W_{i}^{(q)})^{\mT}$, $\bR_{i3} = (R_{i3}^{(1)}, \cdots, R_{i3}^{(q)})^{\mT}$.
Thus, $\Delta\twbX_{2i}^*  = \bW_i + \bR_{i1} + \bR_{i2}  + \bR_{i3} $.
As 
\[\bW_{i} = \sum_{l=1}^{2h-1} \left(1-\frac{|h-l+1|}{h}\right) \Delta \bX_{(2i-2)h + l}^*,
\]
and 
\[
\Delta \bX_{(2i-2)h + l}^* = \begin{pmatrix}
X_{t_{T_{(2i-2)h +l }}}^{(1)} -X_{t_{T_{(2i-2)h +l-1 }}}^{(1)}\\
 \vdots\\
  X_{t_{T_{(2i-2)h +l }}}^{(p)} -X_{t_{T_{(2i-2)h +l-1 }}}^{(p)}\\
  \end{pmatrix},
\]
it reduces to the synchronous setting.
Using a similar argument to that in \eqref{eq: convergence of sum of TWBx}, we have
\[ \lim\limits_{p\to\infty} \frac{\sum_{i=1}^M |\bW_i|^2}{p} = \theta.
\]
Thus, \eqref{eq: xy same require in AT} follows if we can show that there exists a constant $C>0$ such that for large $M$
\begin{equation}\label{eq: min of Wi2}
\min\limits_{1\leq i \leq M}|\bW_{i}|^2\geq C,
\end{equation}
and
\begin{equation}\label{eq: asynchronous key}
\max\limits_{1\leq i \leq M, 1\leq q \leq p}\sqrt{p}|R_{i1}^{(q)} + R_{i2}^{(q)} +  R_{i3}^{(q)}|\to 0,~~~\text{almost surely.}
\end{equation}
Hence, the PA-ATVA matrix based on noisy observations $\Delta\twbY^*_{2i}$ and that based on latent log price $\Delta\twbX^*_{2i}$ have the same LSD.
Further, using \eqref{eq: min of Wi2}, \eqref{eq: asynchronous key}, and Lemma \ref{lem: appendix B lemma 1}, we obtain that
\[
3\frac{\sum_{i=1}^M |\Delta\twbX_{i2}^*|^2}{p} \cdot \frac{p}{\hv}\sum_{i=1}^{\hv}\frac{\Delta\twbX_{i2}^* (\Delta\twbX_{i2}^*)^{\mT}}{|\Delta\twbX_{i2}^*|^2}
~~\text{and}~~3\frac{\sum_{i=1}^M |\bW_i|^2}{p}\cdot\frac{p}{\hv}\sum_{i=1}^{\hv}\frac{\bW_{i}\bW_{i}^{\mT}}{|\bW_{i}|^2}
\]
have the same LSD, and
\[
\lim\limits_{p\to\infty} 3\frac{\sum_{i=1}^M |\Delta\twbX_{i2}^*|^2}{p}=  \lim\limits_{p\to\infty} 3\frac{\sum_{i=1}^M |\bW_i|^2}{p} = \theta.
\]
Moreover, the LSD of  $3 \frac{\sum_{i=1}^M |\bW_{i}|^2}{p}\cdot\frac{p}{\hv}\sum_{i=1}^{\hv}\frac{\bW_{i}\bW_{i}^{\mT}}{|\bW_{i}|^2}$ is generated by Theorem 2.3 of \cite{XiaandZheng20182}.
Therefore, the proof of Theorem \ref{thm: theorem for the TVAPA asynt} is complete.

Finally, we only need to show the proofs of \eqref{eq: min of Wi2} and \eqref{eq: asynchronous key}.
Notice that \eqref{eq: min of Wi2} follows naturally from the proof of \eqref{eq: minimun of twbX}.
To prove \eqref{eq: asynchronous key}, by the normality of $R_{i1}^{(q)}$, $R_{i2}^{(q)}$, and $R_{i3}^{(q)}$, Assumption \eqref{asm:bdd_mu}, and $\max\limits_{1\leq i \leq n, 1\leq q\leq  p} nh(s_{T_i}^{(q)} - t_i) \to 0 $ almost surely, we have
\[
E(R_{i1}^{(q)}) = O(\frac{1}{nh}), ~~~\Var(R_{i1}^{(q)}) = O(\frac{1}{nh^2}),
\]
\[E(R_{i2}^{(q)}) = O(\frac{1}{n}), ~~~\Var(R_{i2}^{(q)}) = O(\frac{1}{nh}),
\]
and
\[E(R_{i3}^{(q)}) = o(\frac{1}{n}), ~~~\Var(R_{i3}^{(q)}) = o(\frac{1}{n}).
\]
Let $R_{i}^{(q)} = R_{i1}^{(q)} + R_{i2}^{(q)} + R_{i3}^{(q)}$; thus, by $C_p$ inequality and the normality of $R_{i}^{(q)} $, we have for any $\kappa\geq 1$, 
\[
 E(R_{i}^{(q)} ) < Cn^{-1} ~~\text{and}~~ E|R_{i}^{(q)} - E(R_{i}^{(q)} ) |^{2\kappa} \leq Cn^{-\kappa},
\]
for all $i, q$.
Hence,  it follows that for any $\varepsilon >0$ and $\kappa \geq 1$,
\[
\begin{split}
&P\left(\max\limits_{i,q} \sqrt{p}| R_{i}^{(q)} - E(R_{i}^{(q)} ) |\geq \varepsilon\right) \\
&\leq \sum\limits_{i,q}\frac{p^{\kappa} E| R_{i}^{(q)} - E(R_{i}^{(q)} ) |^{2\kappa}}{\varepsilon^{2\kappa}}\\
&\leq C\cdot \frac{Mp\cdot p^{\kappa} \cdot n^{-\kappa}}{\varepsilon^{2\kappa}} = O(p^{2-\kappa\beta/(1-\beta)}).
\end{split}
\]
We choose $\kappa$ large enough such that $\kappa\beta/(1-\beta)- 2 > 1$.
Therefore, \eqref{eq: asynchronous key} holds by the Borel--Cantelli lemma.
\proofend

\section{Useful lemmas}
\begin{lem}
(Lemma 1 in \cite{ZhengandLi20112}). Suppose that for each $p$, $\bv_l = (v_l^1, \cdots, v_l^p)^\mathrm{T}$ and $\bm w_l = (w_l^1,\cdots,w_l^p)^{\mathrm{T}},l =1, \cdots, m$, are all $p$-dimensional vectors. Define
\[\bm{\widetilde{S}}_m = \sum_{l=1}^m (\bm v_l + \bm w_l)(\bm v_l + \bm w_l)^\mathrm{T} \text{~and~} \bm S_m = \sum_{l=1}^mw_l(w_l)^\mathrm{T}.\]
Suppose the following conditions are satisfied:
\begin{itemize}
  \item $m = m(p)$ with $\lim_{p\to\infty} p/m = c >0$;
  \item there exists a sequence $\epsilon_p = o(1/\sqrt{p})$ such that for all $p$ and all $l$, all the entries of $\bm v_l$ are bounded by $\epsilon_p$ in absolute value;
  \item $\lim\sup_{p\to\infty}\begin{rm}tr\end{rm}(\bm S_m)/p<\infty$ almost surely.
\end{itemize}
Then $L(F^{\bm{\widetilde{S}}_m}, F^{\bm S_m}) \to 0$ almost surely, where for any two probability distribution functions $F$ and $G$, $L(F, G)$ denotes the Levy distance between them.
\label{lem: appendix B lemma 1}
\end{lem}

\end{supplement}






\bibliographystyle{apa}  


\end{document}